\documentclass[hidelinks,onefignum,onetabnum]{siamart250211}

\usepackage{lipsum}
\usepackage{amsfonts}
\usepackage{graphicx}
\usepackage{epstopdf}
\usepackage{algorithmic}
\ifpdf
\DeclareGraphicsExtensions{.eps,.pdf,.png,.jpg}
\else
\DeclareGraphicsExtensions{.eps}
\fi

\newsiamremark{remark}{Remark}
\newsiamremark{hypothesis}{Hypothesis}
\crefname{hypothesis}{Hypothesis}{Hypotheses}
\newsiamthm{claim}{Claim}
\newsiamremark{fact}{Fact}
\crefname{fact}{Fact}{Facts}
\newsiamremark{example}{Example}

\headers{Linearizing eigenvector nonlinearities}{V. Janssens, K. Meerbergen, and W. Michiels}

\title{Linearizing a nonlinear eigenvalue problem with quadratic rational eigenvector nonlinearities\thanks{This work was supported by the  project G027624N of the Research
		Foundation-Flanders (FWO - Vlaanderen).
}}

\author{Victor Janssens\thanks{Department of Computer Science, KU Leuven, Leuven, Belgium 
		(\email{{victor.janssens@kuleuven.be, karl.meerbergen@kuleuven.be, wim.michiels@cs.kuleuven.be}}).}
	\and Karl Meerbergen\footnotemark[2] \and Wim Michiels\footnotemark[2]}

\usepackage{amsopn}
\usepackage{amssymb}
\usepackage{matlab-prettifier}

\DeclareMathOperator{\rank}{rank}
\DeclareMathOperator{\nrank}{nrank}
\DeclareMathOperator{\nullspace}{null}
\DeclareMathOperator{\col}{\text{col}}

\DeclareMathOperator{\spn}{span}
\DeclareMathOperator{\diag}{diag}
\DeclareMathOperator{\adj}{adj}
\DeclareMathOperator{\trace}{trace}
\DeclareMathOperator{\vect}{\text{vec}}
\DeclareMathOperator{\Vol}{Vol}
\DeclareMathOperator{\Conv}{Conv}

\newcommand{\tran}{\mathsf{T}}
\newcommand{\herm}{\mathsf{H}}

\ifpdf
\hypersetup{
pdftitle={Linearization},
pdfauthor={V. Janssens, K. Meerbergen, and W. Michiels}
}
\fi

\begin{document}

\maketitle

\begin{abstract}
Nonlinear eigenvalue problems with eigenvector nonlinearities (NEPv) are algebraic eigenvalue problems whose matrix depends on the eigenvector. Applications range from computational quantum mechanics to machine learning. Due to its nonlinear behavior, existing methods almost exclusively rely on fixed-point iterations, the global convergence properties of which are only understood in specific cases. Recently, a certain class of NEPv with linear rational eigenvector nonlinearities has been linearized, i.e., the spectrum of the linear eigenvalue problem contains the eigenvalues of the NEPv. This linear problem is solved using structure exploiting algorithms to improve both convergence and reliability. We propose a linearization for a different class of NEPv with quadratic rational nonlinearities, inspired by the discretized Gross-Pitaevskii equation. The eigenvalues of this NEPv form a subset of the spectrum of a linear multiparameter eigenvalue problem which is equivalent to a system of generalized eigenvalue problems expressed in terms of operator determinants. A structure exploiting Arnoldi algorithm is used to filter a large portion of spurious solutions and to accelerate convergence. 
\end{abstract}



\section{Introduction}

We revisit a specific class of the eigenvector dependent nonlinear eigenvalue problem (NEPv), which involves solving
\begin{align}
Av = \lambda B v + \sum_{i = 1}^m f_i(v) C_i v,
\label{eq:1_NEPv_f}
\end{align}
where $A, B, C_i \in \mathbb{C}^{n \times n}$ are complex-valued matrices of size $n$ and where the functions $f_i : \mathbb{C}^n \rightarrow \mathbb{C}$ are scaling invariant for $i \in \{1, 2, 3, \dots, m\}$, i.e., $f(\alpha v) = f(v)$ for all nonzero $v \in \mathbb{C}^n\backslash\{0\}$ and all nonzero $\alpha \in \mathbb{C}\backslash\{0\}$ \cite{CJMU22_linearization}. The eigenpair $(\lambda, v)$ solves this eigenvalue problem if \eqref{eq:1_NEPv_f} holds with $\lambda \in \mathbb{C}$ a complex number and $v \in \mathbb{C}\backslash \{0\}$ a nonzero vector.

More general forms of the NEPv appear in quantum mechanical applications such as electronic structure calculations \cite{KS65_Kohn_Sham, R51_Molecular_Orbital_Theory} and Bose-Einstein condensates \cite{BWM05_GPE_BEC_vortex, HJ24_GPE_NEP}. Other applications in machine learning are p-spectral clustering \cite{BH09_Clustering} and the trace-ratio optimization problem \cite{NBS10_trace_ratio_optimization}. 

Most existing methods rely on fixed point iterations of which the convergence properties are only understood in specific cases. A popular example is the Self-Consistent Field (SCF) iteration, which was originally developed to solve the Hartree and Fock equations in electronic structure calculations \cite{H1927_wave_mechanics_atom_2, R51_Molecular_Orbital_Theory}. The idea is to guess an initial eigenvector and substitute it in the nonlinear part of the NEPv. This gives a generalized eigenvalue problem (GEP) which can be solved to obtain the eigenvector for the next iteration. Repeating this process until the eigenvector does not change any longer after an iteration, gives the solution of this fixed point problem. Convergence of the SCF method is studied, see \cite{CZBL18_NEPv, BLL22_SCF_convergence}, e.g., but it often does not converge in its original form, which is why it is usually combined with acceleration techniques such as Direct Inversion of the Iterative Subspace (DIIS) \cite{P80_DIIS}. Other methods are generalizations of the inverse iteration method \cite{JKM14_Inverse_iteration_NEPv} and contour integration methods \cite{CMT23_Contour_Integration}.

A recent study \cite{CJMU22_linearization} shows that if the functions in the NEPv \eqref{eq:1_NEPv_f} are rational with a linear numerator and denominator in the eigenvector, i.e.
\begin{align}
f_i(v) = \frac{r_i^\tran v}{s_i^\tran v},
\notag
\end{align}
with vectors $r_i, s_i \in \mathbb{C}^n$, the problem can be linearized into a GEP, that is, all eigenvalues of the NEPv are also eigenvalues of the GEP. This is useful as the NEPv can then be solved using efficient and reliable methods developed for the GEP.  In \cite{CJMU22_linearization}, a slight modification to inverse iteration is proposed that exploits the structure of the linearized problem to speed up convergence. 

In this paper, we propose a linearization of \eqref{eq:1_NEPv_f} with one quadratic rational nonlinearity, i.e.,
\begin{align}
\label{eq:1_QF}
f_i(v) = \frac{v^\herm P_i v}{v^\herm Q_i v},
\end{align}
where $m = 1$, $A, C_i, P_i \in H_n = \{X \in \mathbb{C}^{n \times n} | X = X^\herm\}$ are Hermitian and $B, Q_i \in H_n^+ = \{X \in H_n | X \succ 0\}$ positive definite. This type of eigenvector nonlinearity is closely related to a discretization of the Gross-Pitaevskii equation (GPE), also referred to as the nonlinear Schrödinger equation, which needs to be solved in order to obtain the ground state and the wave function of a so called Bose-Einstein condensate (BEC), see for instance \cite{BWM05_GPE_BEC_vortex, HJ24_GPE_NEP} and references therein. In case $C_i$ is a rank one matrix and $C_i = P_i$, the eigenvalues of the NEPv can by obtained by solving an eigenvalue problem with eigenvalue nonlinearities (NEP$\lambda$) instead \cite{JL25_NEPv_to_NEPl} for which efficient methods are available. This approach works very well for this low rank case, whereas our linearization is better suited when the rank of $C_i$ is greater than or equal to $n-1$ and $m = 1$. 

We will show how the NEPv with scalar nonlinearities \eqref{eq:1_QF} for $m = 1$ is related to a multiparameter eigenvalue problem (MEP), which can be constructed directly from the coefficient matrices $A, B, C_1, P_1$ and $Q_1$ and an additional matrix $R$ that can be chosen freely as long as its rank is full. Using operator determinants, the MEP may then be converted into a GEP. The paper proposes a modification to the Arnoldi algorithm \cite{A51_Arnoldi} and exploits the structure of the resulting GEP to obtain the eigenvalues of the NEPv closest to a chosen shift by filtering out some of the spurious solutions.

\medskip 
The paper is organized as follows. An upper bound on the number of eigenvalues of \eqref{eq:1_NEPv_f} with eigenvector nonlinearity \eqref{eq:1_QF} and the link with a system of polynomial equations are discussed in \Cref{sec:problem}. As a first main contribution, \Cref{sec:Linearization} proposes the linearization which contains all the eigenvalues of the NEPv. The resulting problem is solved using structure-exploiting Arnoldi algorithms discussed in \Cref{sec:Numerical_methods}, which is the second contribution. The theoretical results are supported by numerical examples in \Cref{sec:Experiments}, and lastly \Cref{sec:Conclusion} summarizes the results.

\section{Problem properties}
\label{sec:problem}

Given the Hermitian matrices $A, C, P \in H_n$ and positive definite matrices $B, Q \in H_n^+$, find $\lambda \in \mathbb{C}$ and $v \in \mathbb{C}^n\backslash \{0\}$ such that
\begin{align}
Av = \lambda Bv + \frac{v^\herm P v}{v^\herm Q v} Cv.
\label{eq:NEPv_QF}
\end{align}
If such values exist, we say that $(\lambda, v)$ solves the NEPv \eqref{eq:NEPv_QF}. The following lemma states that this NEPv is equivalent to the system of nonlinear equations
\begin{align}
\begin{cases}
M(\lambda, \mu)v = 0, \\
v^\herm S(\mu) v = 0,
\end{cases}
\label{eq:NEPv_mu}
\end{align}
where
\begin{align}
M(\lambda, \mu) & = A - \lambda B - \mu C,
\label{eq:2_M_S}
\\
S(\mu) & = P - \mu Q,
\notag
\end{align}
which we will write in short as $M$ and $S$ if clear from context.

\begin{lemma}
The eigenpair $(\lambda, v)$ solves \eqref{eq:NEPv_QF} if and only if there exists a scalar $\mu \in \mathbb{C}$ such that $(\lambda, \mu, v)$ solves \eqref{eq:NEPv_mu}. Moreover, $\mu$ equals $\frac{v^\herm P v}{v^\herm Q v}$ and both $\lambda$ and $\mu$ are real-valued if they solve the NEPv.
\label{prop:2_NEPv_to_NEPv_mu}
\end{lemma}

\begin{proof}
If $(\lambda, v)$ solves \eqref{eq:NEPv_QF} and we define $\mu = \frac{v^\herm P v}{v^\herm Q v}$, then we have $Av = \lambda Bv + \mu Cv$ and $v^\herm P v = \mu \, v^\herm Qv$ which implies that $(\lambda, \mu, v)$ solves \eqref{eq:NEPv_mu}. 

Conversely, assume $(\lambda, \mu, v)$ solves \eqref{eq:NEPv_mu}, then $v^\herm P v = \mu \, v^\herm Qv$, which can also be written as $\mu = \frac{v^\herm P v}{v^\herm Q v}$ since $v^\herm Q v$ cannot equal zero\footnote{As the eigenvector $v$ is nonzero and $Q$ is positive definite, the value $v^\herm Q v$ is strictly positive.}. Eliminating $\mu$ from the first equation in \eqref{eq:NEPv_mu} gives \eqref{eq:NEPv_QF}.

Because both $P$ and $Q$ are positive definite, the value $\mu = \frac{v^\herm P v}{v^\herm Q v}$ equals its complex conjugate $\overline{\mu}$, and thus $\mu$ is real-valued. Multiplying $\eqref{eq:NEPv_QF}$ with $v^\herm$ on the left and dividing\footnote{This division is allowed due to $v$ being nonzero and $B$ positive definite.} by $v^\herm B v$, gives
\begin{align}
\lambda = \frac{v^\herm (A - \mu C) v}{v^\herm B v},
\notag
\end{align}
and therefore $\lambda$ is real-valued as well because $A - \mu C$ and $B$ are Hermitian.
\end{proof}

In order to obtain an upper bound on the number of isolated solutions of the NEPv \eqref{eq:NEPv_mu}, we show how the eigenvalues are related to the solutions of the system of polynomial equations
\begin{align}
\begin{cases}
\det(M(\lambda, \mu)) = 0, \\
\trace(S(\mu) \, \adj(M(\lambda, \mu))) = 0.
\end{cases}
\label{eq:NEPv_polyn}
\end{align}

\begin{lemma}
If $(\lambda, \mu, v)$ solves the NEPv \eqref{eq:NEPv_mu} with $v$ nonzero, then $(\lambda, \mu)$ solves the equations \eqref{eq:NEPv_polyn}. 
\label{prop:NEPv_polyn}
\end{lemma}

\begin{proof}
Given the solution $(\lambda, \mu, v)$ which satisfies \eqref{eq:NEPv_mu}, we have that $\det(M) = 0$ because $Mv = 0$ for nonzero $v$. To obtain the second equality in \eqref{eq:NEPv_polyn}, we use the fact that the adjugate of the singular matrix $M$ can be written as $\alpha v v^\herm$ for some complex number $\alpha \in \mathbb{C}$ \cite[page 22]{HJ94_Matrix_Analysis}. From \eqref{eq:NEPv_mu}, we have
\begin{align}
v^\herm S v & = 0, \notag \\
\iff \quad \trace(Svv^\herm ) & = 0, \notag \\
\iff \quad \trace(S \, \adj(M)) & = 0,  \notag
\end{align}
and therefore $(\lambda, \mu)$ solves the equations \eqref{eq:NEPv_polyn}.
\end{proof}

The converse is also true under a specific condition and if the eigenvalues are real-valued.

\begin{proposition}
Let $\lambda, \mu \in \mathbb{R}$ be real-valued numbers such that \eqref{eq:NEPv_polyn} holds, and define $\hat{S} = V^\herm S(\lambda, \mu) V$ where the columns of $V \in \mathbb{C}^{n \times k}$ form a basis for the null space of $M(\lambda, \mu)$, then there exists a nonzero vector $v \in \mathbb{C}^n\backslash\{0\}$ such that $(\lambda, \mu, v)$ solves \eqref{eq:NEPv_mu} if and only if $\hat{S}$ is neither positive nor negative definite. Moreover, if $V$ has one column, i.e., $k=1$, then the condition $\pm \hat{S} \not \succ 0$ is automatically satisfied.
\label{prop:polyn_NEPv}
\end{proposition}

\begin{proof}
We start by proving the if and only if statement from left to right. Assume $(\lambda, \mu, v)$ solves the NEPv \eqref{eq:NEPv_mu}, then both equations $Mv = 0$ and $v^\herm S v = 0$ hold. The first equation implies $v \in \nullspace(M)$ and thus there exists a nonzero vector $x \in \mathbb{C}^k\backslash\{0\}$ such that $v = Vx$. The other equation $v^\herm S v = 0$ is equivalent to $x^\herm \hat{S} x = 0$ and thus $\pm \hat{S} \not \succ 0$. 

Conversely, assume $\pm \hat{S} \not \succ 0$ then there exists a nonzero vector $x \in \mathbb{C}^k\backslash\{0\}$ such that $x^\herm \hat{S} x = 0$, which is equivalent to $x^\herm V^\herm S V x = 0$. Note that $MVx = 0$ since $V$ is a basis for the null space of $M$ and therefore $(\lambda, \mu, Vx)$ solves the NEPv \eqref{eq:NEPv_mu}. 

For the case where $k = 1$, there must exist a nonzero $v$ such that $Mv = 0$ and $\adj(M) = vv^\herm$ \cite{HJ94_Matrix_Analysis}. From \eqref{eq:NEPv_polyn} the equation $\trace(S \adj(M)) = 0$ holds and thus $\hat{S} = v^\herm S v = \trace(S v v^\herm) = \trace(S \adj(M)) = 0$ which is the same as $\pm \hat{S} \not \succ 0$.
\end{proof}

\begin{theorem}
\label{prop:2_number_sols}
If the system of polynomial equations \eqref{eq:NEPv_polyn} has a finite number of solutions, the NEPv \eqref{eq:NEPv_QF} has at most $n^2$ eigenvalues.
\end{theorem}

\begin{proof}
From Bézouts theorem, the upper bound on the number of isolated solutions of the system of polynomial equations \eqref{eq:NEPv_polyn} is $n^2$ since both polynomials have a degree of $n$. If $(\lambda, v)$ solves \eqref{eq:NEPv_QF}, then from \cref{prop:2_NEPv_to_NEPv_mu} and \cref{prop:NEPv_polyn}, the pair $(\lambda, \frac{v^\herm P v}{v^\herm Q v})$ solves \eqref{eq:NEPv_polyn}. Therefore an upper bound on the number of eigenvalues of the NEPv is $n^2$.
\end{proof}

An interesting fact is that because $A, B, C, P$ and $Q$ are Hermitian, this system of polynomial equations \eqref{eq:NEPv_polyn} in $\lambda$ and $\mu$ has real-valued coefficients, which can be explained as follows. Let $f(\lambda, \mu) = \det(A - \lambda B - \mu C)$, then for any $\lambda, \mu \in \mathbb{C}$
\begin{align}
\overline{f(\lambda, \mu)} & = \det(\overline{A} - \overline{\lambda} \overline{B} - \overline{\mu} \overline{C}) = \det(A^\tran - \overline{\lambda} B^\tran - \overline{\mu} C^\tran) = f(\overline{\lambda}, \overline{\mu}).
\notag
\end{align} 
As a result, the coefficients of $f$ have an imaginary part equal to zero. Similarly, let $g(\lambda, \mu) = \trace \left((P - \mu Q) \, \adj(A - \lambda B - \mu C)\right)$, then $\overline{g(\lambda, \mu)} = g(\overline{\lambda}, \overline{\mu})$, so the same conclusion can be drawn for the second equation. 

Note that we could solve the NEPv by solving the polynomial system \eqref{eq:NEPv_polyn}, for instance, using homotopy continuation methods \cite{L97_HC}. However, computing its coefficients can be time consuming and numerically unstable because of the evaluation of a determinant and an adjugate in the equations, which further motivates the linearization approach in \Cref{sec:Linearization}.

\begin{example}
\cref{prop:polyn_NEPv} indicates that there may be cases where not all real-valued eigenvalues of the system of polynomial equations \eqref{eq:NEPv_polyn} correspond to eigenvalues of the NEPv \eqref{eq:NEPv_mu}. As an example, consider the matrices
\begin{align}
& A = \begin{bmatrix}
8 & -9 - 6i \\ -9 + 6i & -4
\end{bmatrix}, \quad 
B = \begin{bmatrix}
8 & 3 + 2 i \\ 3 - 2i & 8
\end{bmatrix} , \quad 
C = \begin{bmatrix}
0 & 6 + 4 i \\ 6 - 4 i & 6
\end{bmatrix}
\notag
\\
& P = \begin{bmatrix}
-4 & 6 + 6 i \\ 6 - 6 i & 0
\end{bmatrix}, \quad 
Q = \begin{bmatrix}
6 & -3 - 2 i \\ -3 + 2i & 3
\end{bmatrix}.
\notag
\end{align}
In this case $(\lambda, \mu) = (1, -2)$ is a real-valued solution of the equations
\begin{align}
\begin{cases}
\det(M) = 51 \lambda^2 - 4 \lambda \mu  - 52 \mu^2 - 110 \lambda - 204 \mu - 149 = 0, \\
\trace(S \, \adj(M)) =  98 \lambda \mu + 88 \mu^2 + 92 \lambda + 222 \mu + 196 = 0.
\end{cases}
\notag
\end{align}
Note that $M(1, -2) = 0$, i.e., the dimension of its null space equals 2 and thus the columns of the identity matrix $I$ can be chosen as a basis for the null space of $M$. This means that $\hat{S} = I S(-2) I = P + 2Q$, as defined in \cref{prop:polyn_NEPv}. Since $\hat{S}$ has strictly positive eigenvalues, it is positive definite, which means $\lambda = 1$ is not an eigenvalue of \eqref{eq:NEPv_mu}.
\end{example}

\begin{example}
The upper bound of $n^2$ isolated solutions for the NEPv \eqref{eq:NEPv_QF} can be reached in certain cases. For instance, when $n = 2$ and the NEPv matrices are
\begin{align}
& A = \begin{bmatrix}
4 & 3 + i \\ 3 - i & 1
\end{bmatrix}, \quad 
B = \begin{bmatrix}
16 & 2 - 2 i \\ 2 + 2i & 9
\end{bmatrix} , \quad 
C = \begin{bmatrix}
-8 & 5 - 10 i \\ 5 + 10 i & -17
\end{bmatrix}
\notag
\\
& P = \begin{bmatrix}
6 & -1 + 18 i \\ -1 - 18 i & 4
\end{bmatrix}, \quad 
Q = \begin{bmatrix}
6 & 2 + i \\ 2 - i & 4
\end{bmatrix},
\notag
\end{align}
the problem has 4 different eigenvalues. These are given in \cref{tab:ex_n2_sols} with the corresponding value of $\mu$ and eigenvector $v = \begin{pmatrix}
v_1 & v_2
\end{pmatrix}^\tran$.
\label{ex:2_n2_sols}
\end{example}

\begin{table}[H]
\centering
\begin{tabular}{c||cccc}
\hline
& $\lambda$ & $\mu$ & $v_1$ & $v_2$
\\ \hline
Solution 1 & \phantom{$-$}11.936 & \phantom{$-$}4.0164 & $0.0438 + 0.4424i$ & \phantom{$-$}$0.7073 - 0.5497i$
\\ 
Solution 2 & $-0.0684$ & \phantom{$-$}0.0207 & $0.5209 - 0.0291i$ & $-0.7875 + 0.3282i$
\\ 
Solution 3 & \phantom{$-$}0.1906 & $-1.4229$ & $0.7836 + 0.3013i$ & \phantom{$-$}$0.1508 + 0.5220i$
\\ 
Solution 4 & \phantom{$-$}0.2612 & $-0.3510$ & $0.3963 - 0.7437i$ & \phantom{$-$}$0.4312 - 0.3225i$
\\ \hline
\end{tabular}
\caption{All solution of the NEPv \eqref{eq:NEPv_QF} in \cref{ex:2_n2_sols}.}
\label{tab:ex_n2_sols}
\end{table}

\section{Linearization}
\label{sec:Linearization}

This section makes the link between the NEPv and a GEP. While the linearized problem contains all the eigenvalues (at most $n^2$) of the NEPv, the opposite is not true since the GEP we propose in \Cref{sec:Linearization_Compact} has size $2n^2 - n$, i.e., it may contain eigenvalues that do not correspond to a solution of the NEPv. We will refer to these eigenvalues as spurious solutions and their origin is explained in \Cref{sec:Linearization_Spurious}. Fortunately, these spurious solutions have a special structure in the eigenvectors of the GEP which will be exploited in numerical methods proposed in \Cref{sec:Numerical_methods}.

\subsection{Kronecker linearization}

The NEPv with quadratic eigenvector nonlinearities can be reformulated as a NEPv similar to those studied in \cite{CJMU22_linearization}. This allows us to apply the linearization technique described there to our problem as well. However, the resulting linearized problem is singular and exhibits rapid growth in size with respect to the dimension $n$. While there are workarounds for both issues, the details of which are omitted here, the scaling challenge remains significant. This motivates the more compact linearization approach presented in \Cref{sec:Linearization_Compact}.

Starting from the original NEPv \eqref{eq:NEPv_QF}, the conversion to the different class of NEPv can be realized as follows: Let $\mu = \frac{v^\herm P v}{v^\herm Q v}$ as in \Cref{sec:problem} and define the vector $u \in \mathbb{C}^{n^2}$ as $u = \overline{v} \otimes v$, then $v^\herm P v$ equals $p^\herm u$ where $p = \vect(P) \in \mathbb{C}^{n^2}$ is a vector containing all columns of $P$ stacked on top of each other. Similarly, $v^\herm Q v = q^\herm u$ with $q = \vect(Q)$. From \cref{prop:2_NEPv_to_NEPv_mu}, both $\lambda$ and $\mu$ are real-valued, and therefore
\begin{align}
(\overline{A} \oplus A) u & = (\overline{A} \otimes I + I \otimes A) (\overline{v} \otimes v) = \overline{Av} \otimes v + \overline{v} \otimes Av,
\notag
\\
& = (\lambda \overline{Bv} + \mu \overline{Cv}) \otimes v + \overline{v} \otimes (\lambda Bv + \mu Cv),
\notag
\\
& = \lambda (\overline{C} \oplus C) u + \mu (\overline{C} \oplus C) u.
\notag
\end{align}
Eliminating $\mu$ from the equations gives the NEPv with eigenvector nonlinearities described in \cite{CJMU22_linearization}:
\begin{align}
(\overline{A} \oplus A) u = \lambda (\overline{C} \oplus C) u + \frac{p^\herm u}{q^\herm u} (\overline{C} \oplus C) u.
\notag
\end{align}
Linearizing this NEPv results in a GEP of size $n^4$ which scales worse than the linearization of size $2n^2 - n$ proposed in the next section. Therefore, the remaining parts of this paper will focus on the latter.

\subsection{Compact linearization}
\label{sec:Linearization_Compact}

First, the NEPv \eqref{eq:NEPv_QF} is converted into a MEP. Let $R \in \mathbb{C}^{n \times (n-1)}$ be a matrix of full column rank, then the linearization is defined as the MEP
\begin{align}
\begin{cases}
M(\lambda, \mu)v = 0, \\
\begin{bmatrix}
0 & R^\herm M(\lambda, \mu) \\ M(\lambda, \mu)R & S(\mu)
\end{bmatrix}w = 0.
\end{cases}
\label{eq:3_MEP_R}
\end{align}
This is a two-parameter eigenvalue problem with solution $(\lambda, \mu, v, w)$, i.e., eigentuple $(\lambda, \mu) \in \mathbb{C} \times \mathbb{C}$ and eigenvectors $v \in \mathbb{C}^n\backslash\{0\}$, $w \in \mathbb{C}^{2n-1}\backslash\{0\}$, which is clearer to see if we rewrite \eqref{eq:3_MEP_R} using $M(\lambda, \mu) = A - \lambda B - \mu C$ and $S(\mu) = P - \mu Q$ as
\begin{align}
\begin{cases}
Av = \lambda Bv + \mu Cv, \\
\hat{A}w = \lambda \hat{B}w + \mu \hat{C}w, \\
\end{cases}
\notag
\end{align}
where
\begin{align}
\hat{A} = \begin{bmatrix}
0 & R^\herm A \\ AR & P
\end{bmatrix}, \quad \hat{B} = \begin{bmatrix}
0 & R^\herm B \\ BR & 0
\end{bmatrix}, \quad \hat{C} = \begin{bmatrix}
0 & R^\herm C \\ CR & Q
\end{bmatrix}.
\notag
\end{align}
In general, the number of isolated solutions equals $2n^2 - n$ \cite{A72_MEP}. In what follows, we will show how the NEPv \eqref{eq:NEPv_mu} is related to this MEP.

\begin{theorem}
\label{prop:NEPv2MEP}
If $(\lambda, v)$ solves the NEPv \eqref{eq:NEPv_QF}, then there exists a scalar $\mu \in \mathbb{R}$ and a vector $w \in \mathbb{C}^{2n-1} \backslash \{0\}$ such that $(\lambda, \mu, v, w)$ solves the MEP \eqref{eq:3_MEP_R} and such that $w$ equals $\begin{bmatrix}
w_1^\tran & \alpha v^\tran
\end{bmatrix}^\tran$ for some $w_1 \in \mathbb{C}^{n-1}$ and $\alpha \in \mathbb{C}$.
\end{theorem}

\begin{proof}
Suppose $(\lambda, v)$ solves \eqref{eq:NEPv_QF}, then it follows from \cref{prop:2_NEPv_to_NEPv_mu} that $\lambda$ is real-valued and that there exists a scalar $\mu \in \mathbb{R}$ such that
\begin{align}
\begin{cases}
M(\lambda, \mu)v = 0, \\
v^\herm S(\mu) v = 0,
\end{cases}
\notag
\end{align}
where $M$ and $S$ are defined in \eqref{eq:2_M_S}. Because $M$ is Hermitian, the vector $v$ must lie in the left null space of the square matrix $\begin{bmatrix}
MR & Sv
\end{bmatrix}$, i.e., $v^\herm \begin{bmatrix}
MR & Sv
\end{bmatrix} = 0$ with $R \in \mathbb{C}^{n \times (n-1)}$ a matrix of full rank. As a result, there exists a vector $ w_1 \in \mathbb{C}^{n-1}$ and a scalar $\alpha \in \mathbb{C}$ such that $\begin{bmatrix}
w_1^\tran & \alpha
\end{bmatrix}^\tran$ lies in the right null space of this matrix, that is, $MR w_1 + \alpha S v = 0$. Also since $Mv = 0$, the equality $\alpha R^\herm Mv = 0$ holds and therefore
\begin{align}
\begin{cases}
Mv = 0, \\
\begin{bmatrix}
0 & R^\herm M \\ MR & S
\end{bmatrix}\begin{bmatrix}
w_1 \\ \alpha v
\end{bmatrix} = 0,
\end{cases}
\notag
\end{align}
meaning $(\lambda, \mu, v, \begin{bmatrix}
w_1^\tran & \alpha v^\tran
\end{bmatrix}^\tran)$ solves the MEP \eqref{eq:3_MEP_R}.
\end{proof}

We thus have a linear problem of which the solutions contain all eigenvalues of the original nonlinear problem. Of course, since our MEP has more isolated solutions than the NEPv, the converse cannot be true unless some assumptions are specified. 

\begin{proposition}
Assume $\lambda $ and $\mu$ are real-valued, that $\rank(M(\lambda, \mu))$ equals $\rank(M(\lambda, \mu)R)$ and that $w = \begin{bmatrix}
w_1^\tran & w_2^\tran
\end{bmatrix}^\tran \in \mathbb{C}^{2n-1}$ is some vector such that $w_2 \in \mathbb{C}^{n}$ is nonzero. If $(\lambda, \mu, v, w)$ solves the MEP \eqref{eq:3_MEP_R}, then $(\lambda, w_2)$ solves \eqref{eq:NEPv_QF}.
\label{prop:3_MEP_to_NEPv}
\end{proposition}

\begin{proof}
Since $(\lambda, \mu, v, w)$ solves \eqref{eq:3_MEP_R} the equations
\begin{subequations}
\label{eq:3_MEP_expanded}
\begin{align}
& Mv = 0, \\
& R^\herm M w_2 = 0, 
\label{eq:3_MEP_expanded_b}
\\
& MR w_1 + S w_2 = 0,
\label{eq:3_MEP_expanded_c}
\end{align}
\end{subequations}
must hold. From the assumptions we have that $M$ is Hermitian and that the ranks of $M$, $MR$ and $R^\herm M$ are equal. By consequence, $M$ and $R^\herm M$ must share the same null space, hence $M w_2 = 0$ and $w_2^\herm M = 0$ because of \eqref{eq:3_MEP_expanded_b}. The vector $w_2$ cannot be the zero vector and therefore $(\lambda, \mu, w_2)$ solves
\begin{align}
\notag
\begin{cases}
M w_2 = 0, \\
w_2^\herm S w_2 = 0.
\end{cases}
\end{align}
\cref{prop:2_NEPv_to_NEPv_mu} then states that $(\lambda, w_2)$ is a solution of \eqref{eq:NEPv_QF}.
\end{proof}

\begin{remark}
Note that the assumption that $w_2 \ne 0$ in the last proposition is automatically satisfied if $\rank(MR) = n-1$. Indeed, if we assume $w_2$ is zero anyway, we have that $MRw_1 = 0$ from \eqref{eq:3_MEP_expanded_c} which is only possible if $w_1$ is zero since $MR$ is of full rank. But now the vector $w$ is zero which contradicts with the fact that eigenvectors must be nonzero and we must deduce that $w_2 \ne 0$. 
\end{remark}

\begin{remark}
If $(\lambda, \mu, v, w)$ solves the MEP \eqref{eq:3_MEP_R}, then the matrix $M(\lambda, \mu)$ is singular. Hence, the assumption that the rank of $M$ must equal the rank of $MR$ in \cref{prop:3_MEP_to_NEPv} is satisfied with probability one for a randomly chosen matrix $R \in \mathbb{C}^{n \times (n-1)}$. 
\end{remark}

The linearization \eqref{eq:3_MEP_R} discussed up until now can thus be solved using MEP solving techniques \cite{HKP01_Jacobi_Davidson_2_param_MEP, RMM22_subspace_for_MEP_tensor_based, DYY16_HC_MEP, RJ21}, but in this article we will convert it to a GEP that is solved in \Cref{sec:Numerical_methods} to filter out a portion of the spurious solutions. This conversion can be obtained after defining the operator determinants
\begin{align}
\begin{cases}
\Delta_0 = B \otimes \hat{C} - C \otimes \hat{B}, \\
\Delta_1 = A \otimes \hat{C} - C \otimes \hat{A}, \\
\Delta_2 = B \otimes \hat{A} - A \otimes \hat{B}. 
\end{cases}
\notag
\end{align}

\begin{theorem}
\label{thm:MEP2GEP}
\cite[Chapter 6]{A72_MEP} If $\Delta_0$ is nonsingular, the MEP \eqref{eq:3_MEP_R} is equivalent to the system of generalized eigenvalue problems
\begin{align}
\label{eq:3_GEP_system}
\begin{cases}
\Delta_1 z = \lambda \Delta_0 z, \\
\Delta_2 z = \mu \Delta_0 z,
\end{cases}
\end{align}
with $z = v \otimes w$.
\end{theorem}
Since we are interested in the eigenvalue $\lambda$, it suffices to solve the first GEP in \eqref{eq:3_GEP_system} 
\begin{align}
\Delta_1 z = \lambda \Delta_0 z,
\label{eq:3_GEP}
\end{align}
which will be referred to as the linearization of the class of NEPv \eqref{eq:NEPv_QF}. The methods developed in \Cref{sec:Numerical_methods} solve this GEP using the Arnoldi iteration to obtain the NEPv eigenvalues. In order for this to work, the pencil must be nonsingular which is true if $\Delta_0$ is nonsingular as stated in \cref{thm:MEP2GEP}. The following lemma and theorem shed light on the invertibility of $\Delta_0$, but before providing both statements, recall that
\begin{align}
\hat{B} = \begin{bmatrix}
0 & R^\herm B \\ BR & 0
\end{bmatrix} \quad \text{and} \quad \hat{C} = \begin{bmatrix}
0 & R^\herm C \\ CR & Q
\end{bmatrix}.
\notag
\end{align}

\begin{lemma}
\label{lemma:3_ChBh_pencil}
Let $\lambda \in \mathbb{C}$ be some complex number. There exists a nonzero vector $z \in \mathbb{C}^{2n-1}\backslash \{0\}$ such that $(\lambda, z)$ solves $\hat{C}z = \lambda \hat{B}z$ with $\lambda \in \mathbb{R}$ if and only if $(\lambda, x)$ solves $CRx = \lambda BRx$ for some $x \in \mathbb{C}^{n-1}\backslash\{0\}$.
\end{lemma}

\begin{proof}
We start by proving the statement from right to left: Assume $(\lambda, x)$ solves $CRx = \lambda BRx$, then $(\lambda, Rx)$ is an eigenpair of the pencil $(C, B)$ as $R$ is a full rank matrix. Since both $B$ and $C$ are Hermitian and $B$ is positive definite, $\lambda$ must be real-valued. Define $z = \begin{bmatrix}
x^\tran & 0^\tran
\end{bmatrix}^\tran \in \mathbb{C}^{2n-1}$, then $\hat{C}z = \lambda \hat{B}z$ and thus $\lambda$ is a real-valued eigenvalue of the pencil $(\hat{C}, \hat{B})$.

The converse is also true: Assume $(\lambda, z)$ solves $\hat{C}z = \lambda \hat{B}z$ with $\lambda$ real-valued and $z = \begin{bmatrix}
x^\tran & y^\tran
\end{bmatrix}^\tran$ nonzero, then the following two equations must hold:
\begin{align}
\begin{cases}
R^\herm Cy = \lambda R^\herm B y, \\
CRx + Qy = \lambda BRx.
\end{cases}
\notag
\end{align}
Because $Q \succ 0$, it has a Cholesky factorization $Q = LL^\herm$. Left multiplying the first equation with $x^\herm$ and defining $b = L^{-1} BRx$ and $c = L^{-1} CRx$ gives
\begin{align}
\begin{cases}
c^\herm L^\herm y = \lambda b^\herm L^\herm y, \\
L^\herm y = \lambda b - c.
\end{cases}
\label{eq:3_ChBh_pencil_after_cholesky}
\end{align}
Eliminating $L^\herm y$ and rearranging terms yields
\begin{align}
\|b\|^2 \lambda^2 - 2 \Re(b^\herm c) \lambda + \|c\|^2 = 0,
\label{eq:3_quadratic_eq}
\end{align}
where $\Re(b^\herm c)$ is the real part of $b^\herm c$. This last equation is always quadratic in $\lambda$ as $b$ cannot be the zero vector due to the following reasoning: Assume $b = 0$ then from the definition of $b$, $x = R^\dagger B^{-1}Lb = 0$ because $R \in \mathbb{C}^{n \times (n-1)}$ is of full rank, and $c = 0$ from equation \eqref{eq:3_quadratic_eq}. But now $y = 0$ because of equation \eqref{eq:3_ChBh_pencil_after_cholesky} and thus $z = 0$ which is not possible because this is a (nonzero) eigenvector of $\hat{C} - \lambda \hat{B}$, thus $b \ne 0$. Solving equation \eqref{eq:3_quadratic_eq} for $\lambda$ gives
\begin{align}
\lambda = \frac{\Re(b^\herm c) \pm \sqrt{D}}{\|b\|^2},
\notag
\end{align}
where $D = \Re(b^\herm c)^2 - \|b\|^2 \|c\|^2$ is the discriminant. Using the Cauchy-Schwarz inequality, we obtain
\begin{align}
\Re (b^\herm c)^2 \leq |b^\herm c|^2 \leq \|b\|^2\|c\|^2,
\end{align}
and therefore $D \leq 0$ meaning the only way to get a real-valued eigenvalue of this problem is by setting the discriminant equal to zero. This is true if and only if $c = \alpha b$ with $\alpha \in \mathbb{R}$ \cite[Theorem 5.1.4]{HJ94_Matrix_Analysis}, in which case $\lambda = \alpha$. Using our definitions of $b$ and $c$, we can replace $c = \lambda b$ with
\begin{align}
L^{-1}CRx = \lambda L^{-1}BRx,
\notag
\end{align}
and thus $(\lambda, x)$ solves $CRx = \lambda BRx$.
\end{proof}

A direct consequence of this lemma is that if the rectangular pencil $(CR, BR)$ has no solutions, all eigenvalues of $(\hat{C}, \hat{B})$ have a nonzero imaginary part. This fact is indirectly used in the following theorem.

\begin{theorem}
The operator determinant $\Delta_0$ is singular if and only if there exists a value $\lambda \in \mathbb{C}$ and a nonzero vector $x \in \mathbb{C}^{n-1}\backslash \{0\}$ such that $CRx = \lambda BRx$. 
\label{thm:singular_Delta_0}
\end{theorem}

\begin{proof}
We start the proof from right to left: Assume $CRx = \lambda Bx$ is solved for some $\lambda \in \mathbb{C}$ and nonzero $x \in \mathbb{C}^{n-1}\backslash\{0\}$, then \cref{lemma:3_ChBh_pencil} tells us that $\lambda$ must be real-valued and that $\hat{C}y = \lambda \hat{B}y$ for some nonzero vector $y \in \mathbb{C}^{2n-1}\backslash \{0\}$. Define the vector $z = Rx \otimes y$, then
\begin{align}
\Delta_0 z & = BRx \otimes \hat{C} y - CRx \otimes \hat{B} y,
\notag
\\
& = \lambda BRx \otimes \hat{B} y - \lambda BRx \otimes \hat{B} y = 0,
\notag
\end{align}
and since $z$ is nonzero, $\Delta_0$ must be singular.

Conversely, if $\Delta_0$ is singular, there exists a nonzero vector $z \in \mathbb{C}^{2n^2 - n}\backslash \{0\}$ such that $\Delta_0 z = 0$. Define $Z \in \mathbb{C}^{(2n-1) \times n}$ such that $z = \vect(Z)$, then $\Delta_0 z = 0$ can be written as a Sylvester equation in $Z$:
\begin{align}
\hat{C}ZB^\tran - \hat{B}ZC^\tran = 0
\notag
\end{align}
Such an equation has a nonzero solution if and only if the pencils $(C, B)$ and $(\hat{C}, \hat{B})$ have a common eigenvalue\footnote{The correct statement says there is a unique solution if and only if the two pencils have disjoint spectra and they are both regular. However, in our case $C - \lambda B$ is regular since $B$ is invertible, and $\hat{C} - \lambda \hat{B}$ is regular because otherwise \cref{lemma:3_ChBh_pencil} states that $C - \lambda B$ is singular which is impossible.} \cite{GLAM92_general_sylvester_method}, that is, there exists a $\lambda$ such that both $C - \lambda B$ and $\hat{C} - \lambda \hat{B}$ are rank deficient. This eigenvalue must be real-valued due to $C$ being Hermitian and $B$ positive definite, and consequently, there exists a nonzero vector $x \in \mathbb{C}^{n-1}\backslash\{0\}$ such that $CRx = \lambda BRx$ because of \cref{lemma:3_ChBh_pencil}.
\end{proof}

A problem occurs when $\rank(C) < n - 1$, in which case $CR$ is rank deficient meaning there must be some nonzero vector $x \in \mathbb{C}^{n-1}$ such that $CRx = 0$. More precisely, $(\lambda = 0, x)$ solves the equation $CRx = \lambda BRx$ and $\Delta_0$ must be singular because of \cref{thm:singular_Delta_0}, so the assumption of \cref{thm:MEP2GEP} does not hold. Even worse, the vector $Rx \otimes \begin{bmatrix}
x^\tran & 0
\end{bmatrix}^\tran \in \mathbb{C}^{n(2n-1)}$ lives in the null space of both $\Delta_0$ and $\Delta_1$, so the entire pencil $(\Delta_1, \Delta_0)$ is singular. This case is further investigated in \Cref{sec:Numerical_methods_singular} as we are still able to find the NEPv eigenvalues using this linearization \eqref{eq:3_GEP}.

\subsection{Spurious solutions}
\label{sec:Linearization_Spurious}

The MEP \eqref{eq:3_MEP_R} and its corresponding GEP \eqref{eq:3_GEP_system} have at most $2n^2-n$ isolated solutions, which exceeds the upper bound of $n^2$ isolated solutions of the NEPv \eqref{eq:NEPv_QF}, so after solving the linearized problem we should expect some false eigenvalues to appear. Part of the solutions change if a different choice for $R$ is taken in the linearization, and because of their dependence on $R$, they must be spurious. More specifically, consider the problem of finding the numbers $\lambda, \mu \in \mathbb{C}$ and the nonzero vector $x \in \mathbb{C}^{n-1}\backslash\{0\}$ such that
\begin{align}
\label{eq:3_rMEP1}
M(\lambda, \mu)Rx = 0,
\end{align}
i.e., we are looking for a scalar $\mu$ such that there is an eigenvalue of $(A - \mu C, B)$ with an eigenvector in the column space of $R$. If we use the fact that $M = A - \lambda B - \mu C$, it becomes clear that we are essentially trying to solve the rectangular multiparameter eigenvalue problem (rMEP)
\begin{align}
\notag
ARx = \lambda BRx + \mu CRx.
\end{align}
In the generic case, this problem has exactly $\ell = \frac{1}{2}n(n-1)$ eigenvalues counting multiplicities \cite[Lemma 1]{HKP23_rMEP, SS09_linear_rMEP}. Note that because $(\lambda, \mu, x)$ solves the rMEP, the equation $MRx = 0$ holds, which also means that
\begin{align}
\notag
\begin{cases}
M(Rx) = 0, \\
\begin{bmatrix}
0 & R^\herm M \\
MR & S
\end{bmatrix}
\begin{bmatrix}
x \\ 0
\end{bmatrix} = 0.
\end{cases}
\end{align}
The vector $Rx$ cannot be zero since $R$ is of full rank and thus $(\lambda, \mu, Rx, \begin{bmatrix}
x^\tran & 0
\end{bmatrix}^\tran)$ solves our linearization \eqref{eq:3_MEP_R}.

A similar reasoning can be said about the rMEP,
\begin{align}
\notag
y^\herm R^\herm M(\lambda, \mu) = 0
\end{align}
which also has exactly $\ell$ solutions $(\lambda, \mu, y)$ counting multiplicities in the generic case. Indeed,
\begin{align}
\notag
\begin{cases}
(Ry)^\herm M = 0 \\
\begin{bmatrix}
y^\herm & 0
\end{bmatrix}
\begin{bmatrix}
0 & R^\herm M \\
MR & S
\end{bmatrix} = 0
\end{cases}
\end{align}
and therefore $(\lambda, \mu)$ is also an eigenvalue of \eqref{eq:3_MEP_R}.

The two rMEPs discussed above explain the appearance of spurious eigenvalues that depend on $R$ when solving $\Delta_1 z = \lambda \Delta_0 z$ and in \Cref{sec:Numerical_methods_nonsingular} we use a property about their right and left eigenvectors in order to filter them out. However, these eigenvalues may not account for all the spurious solutions. Recall that the upper bound $n^2$ on the number of isolated solutions of $\eqref{eq:NEPv_QF}$ is based on the number of solutions of the polynomial equations \eqref{eq:NEPv_polyn}. It turns out that all solutions of these equations that do not correspond to eigenvalues of the NEPv \eqref{eq:NEPv_QF}, such as the complex conjugate pairs, are also spurious solutions, which is a consequence of \cref{prop:polyn_2_MEP}.

\begin{proposition}
If $(\lambda, \mu)$ solves the polynomial equations \eqref{eq:NEPv_polyn}, then there exist nonzero vectors $v \in \mathbb{C}^n$ and $w \in \mathbb{C}^{2n-1}$ such that $(\lambda, \mu, v, w)$ solves the MEP \eqref{eq:3_MEP_R}.
\label{prop:polyn_2_MEP}
\end{proposition}

\begin{proof}
Given $(\lambda, \mu)$ such that \eqref{eq:NEPv_polyn} holds, then $k = \dim(\nullspace(M)) \geq 1$ since $\det(M) = 0$. We distinguish two cases:

\begin{itemize}
\item If $k = 1$, there exist nonzero vectors $u, v \in \mathbb{C}^n \backslash \{0\}$ such that $Mv = 0$, $u^\herm M = 0$ and $\adj(M) = vu^\herm$ \cite[page 22]{HJ94_Matrix_Analysis}. Using \eqref{eq:NEPv_polyn} we then have $0 = \trace(S \, \adj(M)) = \trace(S v u^\herm) = u^\herm S v$.

\item In case $k > 1$, let the matrices $X \in \mathbb{C}^{n \times k}$ and $Y \in \mathbb{C}^{n \times k}$ be such that their columns form a basis for the left and right null space of $M$, respectively, and define $\hat{S} = X^\herm S Y \in \mathbb{C}^{k \times k}$. Choose any nonzero $x \in \mathbb{C}^k\backslash\{0\}$ and $y \in \mathbb{C}^k\backslash\{0\}$ such that $x^\herm \hat{S} y = 0$ which is always possible since $k > 1$, then $Mv = 0$, $u^\herm M = 0$ and $u^\herm S v = 0$ with $u = Xx$ and $v = Yy$.
\end{itemize}

From the equalities $u^\herm M = 0$ and $u^\herm  S v = 0$ we have that $u$ lies in the left null space of the square matrix $\begin{bmatrix}
MR & Sv
\end{bmatrix}$. Let $\begin{bmatrix}
w_1^\tran & \alpha
\end{bmatrix}^\tran \in \mathbb{C}^n$ be a nonzero vector in its right null space, then it is easy to check that $(\lambda, \mu, v, \begin{bmatrix}
w_1^\tran & \alpha v^\tran
\end{bmatrix}^\tran)$ solves \eqref{eq:3_MEP_R}.
\end{proof}

\section{Numerical methods}
\label{sec:Numerical_methods}

Before proposing numerical methods to solve the resulting linearized GEP, we must know whether our pencil $(\Delta_1, \Delta_0)$ is singular. In case it is nonsingular, we may rely on well developed algorithms in the literature such as the Arnoldi algorithm which will be adapted in two ways such that it filters a large portion of the spurious solutions in \Cref{sec:Numerical_methods_nonsingular}. In case the pencil is singular, the problem becomes more difficult to solve. \Cref{sec:Numerical_methods_singular} studies the case where $\rank(C) < n-1$ as it induces such a singular pencil. For these problems, the relation between the MEP and the system of GEPs is less understood \cite{KP22_singular_2_MEP_2}, but the eigenvalues can still be retrieved using projection, augmentation or rank-completing perturbations \cite{HMP23_Singular_GEP_part_II, MW24_shift_and_invert_Arnoldi} and for this specific problem we can show how it may be solved using the methods from \Cref{sec:Numerical_methods_nonsingular} as well.

\subsection{Nonsingular pencil}
\label{sec:Numerical_methods_nonsingular}

To obtain the eigenvalues of the NEPv \eqref{eq:NEPv_QF}, it suffices to select the real-valued solutions of the generalized eigenvalue problem
\begin{align}
\notag
\Delta_1 z = \lambda \Delta_0 z.
\end{align}
Finding all $2n^2 - n$ eigenvalues is possible, for instance by using the QZ algorithm, but this is applicable to small-scale problems only, and usually we are interested in one or a couple of eigenvalues close to a shift for which faster algorithms are available. This section focuses on Krylov methods, specifically the Arnoldi algorithm which will also search for the undesirable spurious solutions making it expensive for large problems. To elevate this problem, we choose a starting vector for the algorithm such that all iterates are restricted to an invariant subspace which excludes a portion of the spurious eigenvalues.

\subsubsection{Filtering Arnoldi method}

Recall that the Arnoldi algorithm \cite{A51_Arnoldi} builds and orthogonalizes the Krylov space
\begin{align}
\notag
\mathcal{K}_k = \{z_0, (\Delta_1 - \sigma \Delta_0)^{-1} \Delta_0 z_0, \dots, ((\Delta_1 - \sigma \Delta_0)^{-1} \Delta_0)^{k-1} z_0 \},
\end{align}
where $z_0$ is the initial starting vector and $\sigma$ a shift. The orthonormal basis is given by $Z_k \in \mathbb{C}^{(2n^2 - n) \times k}$ such that $\mathcal{K}_k = \spn(Z_k)$. The Ritz values are the eigenvalues of the projected problem in Hessenberg form
\begin{align}
H_k = Z_k^\herm (\Delta_1 - \sigma \Delta_0)^{-1} \Delta_0 Z_k.
\label{eq:4_Hessenberg_proj}
\end{align}
This method can be improved by observing the following fact: \cref{prop:NEPv2MEP} and \cref{thm:MEP2GEP} show that an eigenvalue of the NEPv \eqref{eq:NEPv_mu} corresponds to an eigenvector $z = v \otimes \begin{bmatrix}
w_1^\tran & \alpha v^\tran
\end{bmatrix}^\tran$ of the linear GEP \eqref{eq:3_GEP_system}. This vector has a special structure and is an element of the vector space
\begin{align}
\label{eq:3_vector_space}
\mathcal{Z} = \left\{z \in \mathbb{C}^{(2n^2 - n)} \, \Big| \, z = \vect \left( \begin{bmatrix}
W \\ V
\end{bmatrix} \right),\, V = V^\tran \in \mathbb{C}^{n \times n},\, W \in \mathbb{C}^{(n-1) \times n} \right\}
\end{align}
of dimension $n^2 + \ell$ where $\ell = \frac{1}{2}n(n-1)$. If the Arnoldi algorithm is applied by choosing the starting vector $z_0 \in \mathcal{Z}$, then $\ell$ of the spurious eigenvalues can be filtered out. We prove this in \cref{col:4_invariant_iteration} using \cref{lemma:4_inv_subspace} and the two technical lemmas in \Cref{sec:app_tech_lemmas}.

\begin{lemma}
Assume the matrices $\Delta_0$, $\Delta_1 - \sigma \Delta_0$ and $\mathbf{A} - \sigma \mathbf{B} \in \mathbb{C}^{(\frac{n^2 - n}{2}) \times (\frac{n^2 - n}{2})}$ are all nonsingular with
\begin{align*}
\begin{cases}
\mathbf{A} = L^\tran (R^\herm \otimes R^\herm) (A \otimes C - C \otimes A)T, \\
\mathbf{B} = L^\tran (R^\herm \otimes R^\herm) (B \otimes C - C \otimes B)T,
\end{cases} 
\end{align*}
where $L \in \mathbb{C}^{(n-1)^2 \times (\frac{n^2 - n}{2})}$ and $T\in \mathbb{C}^{n^2 \times (\frac{n^2 - n}{2})}$ are full rank matrices defined as
\begin{align*}
\begin{cases}
\col(L) = \{ e_j^{(n-1)} \otimes e_i^{(n-1)} \mid 1 \leq i \leq j \leq n-1 \}, \\
\col(T) = \{ e_j^{(n)} \otimes e_i^{(n)} - e_i^{(n)} \otimes e_j^{(n)} \mid 1 \leq i < j \leq n \},
\end{cases}
\end{align*}
and $e_i^{(m)} \in \mathbb{R}^m$ is the $i$-th standard unit vector. If $z_{k} \in \mathcal{Z}$ and $z_{k+1} = (\Delta_1 - \sigma \Delta_0)^{-1}\Delta_0 z_{k}$, then $z_{k+1} \in \mathcal{Z}$.
\label{lemma:4_inv_subspace}
\end{lemma}

\begin{proof}
First, we rewrite the equation $(\Delta_1 - \sigma \Delta_0)z_{k+1} = \Delta_0 z_{k}$ as
\begin{align}
\hat{C} Z_{k+1} \Gamma^\tran - \hat{\Gamma}Z_{k+1} C^\tran = \hat{C} Z_{k} B^\tran - \hat{B}Z_{k} C^\tran,
\label{eq:4_system_for_invariant_subspace_proof}
\end{align}
where $\Gamma = A - \sigma B$, $\hat{\Gamma} = \hat{A} - \sigma \hat{B}$, $z_{k} = \vect(Z_{k})$ and $z_{k+1} = \vect(Z_{k+1})$. Let 
\begin{align*}
Z_{k} = \begin{bmatrix}
W_{k} \\ V_{k}
\end{bmatrix} \quad \text{and} \quad Z_{k+1} = \begin{bmatrix}
W_{k+1} \\ V_{k+1}
\end{bmatrix},
\end{align*}
then the first $n-1$ rows of equation \eqref{eq:4_system_for_invariant_subspace_proof} read
\begin{align*}
R^\herm ( CV_{k+1} \Gamma^\tran - \Gamma V_{k+1} C^\tran) = R^\herm ( C V_{k} B^\tran - B V_{k} C^\tran).
\end{align*}
Right multiplying with $\overline{R}$ and defining $F = R^\herm ( C V_{k} B^\tran - B V_{k} C^\tran) \overline{R}$ yields
\begin{align*}
R^\herm ( CV_{k+1} \Gamma^\tran - \Gamma V_{k+1} C^\tran) \overline{R} = F.
\end{align*}
Since $z_{k} \in \mathcal{Z}$, $V_{k}$ is a symmetric matrix and as a result $F$ is skew symmetric, i.e., the symmetric matrix $F + F^\tran$ equals zero and therefore
\begin{align*}
R^\herm ( C (V_{k+1} - V_{k+1}^\tran) \Gamma^\tran - \Gamma (V_{k+1} - V_{k+1}^\tran) C^\tran) \overline{R} = F + F^\tran = 0.
\end{align*}
This is a homogeneous system of linear equations in $X = V_{k+1} - V_{k+1}^\tran$ that may equivalently be written in terms of Kronecker products as
\begin{align*}
(R^\herm \otimes R^\herm)(\Gamma \otimes C - C \otimes \Gamma)\vect(X) = 0.
\end{align*}
Due to the symmetry, the equation has redundant rows that can be eliminated by left multiplying with $L^\tran$. Moreover, $X$ is skew symmetric and can be written as $\vect(X) = T x$ where $x \in \mathbb{C}^{n(n-1)/2}$ contains the strictly upper triangular elements of $X$. Combining these reductions gives
\begin{align*}
L^\tran (R^\herm \otimes R^\herm)(\Gamma \otimes C - C \otimes \Gamma)Tx = 0,
\end{align*}
and therefore
\begin{align*}
(\mathbf{A} - \sigma \mathbf{B})x = 0.
\end{align*}
We made the assumption that this matrix is invertible, that is, the unique solution is the zero vector $x = 0$ which implies $X = V_{k+1} - V_{k+1}^\tran = 0$. By consequence, $V_{k+1}$ is symmetric and $z_{k+1} \in \mathcal{Z}$.
\end{proof}

\begin{theorem}
Assume $\Delta_0$ is nonsingular and that $\sigma \in \mathbb{C}$ is not an eigenvalue of the pencil $\Delta_1 - \lambda \Delta_0$. Let $z_{k} \in \mathcal{Z}$ and $z_{k+1} = (\Delta_1 - \sigma \Delta_0)^{-1}\Delta_0 z_{k}$. If either 
\begin{itemize}
\item (1) the square pencil $(\mathbf{A}, \mathbf{B})$ of size $\frac{1}{2}n(n-1)$ is nonsingular with $\mathbf{A}$ and $\mathbf{B}$ defined in \cref{lemma:4_inv_subspace}, or
\item (2) the system of GEPs \eqref{eq:3_GEP_system} has $2n^2 - n$ distinct eigentuples $(\lambda, \mu)$,
\end{itemize}
then $z_{k+1} \in \mathcal{Z}$.
\label{col:4_invariant_iteration}
\end{theorem}

\begin{proof}
We adopt two complementary approaches. Our first approach is based on directly analyzing the matrix vector product for $z_{k+1}$. Assume that condition (1) holds and choose a $\sigma_0$, not necessarily equal to the shift $\sigma$, such that both $\Delta_1 - \sigma_0 \Delta_0$ and $\mathbf{A} - \sigma_0 \mathbf{B}$ are nonsingular, then from \cref{lemma:4_inv_subspace} the subspace $\mathcal{Z}$ is invariant under multiplication with the matrix $\Delta = (\Delta_1 - \sigma_0 \Delta_0)^{-1} \Delta_0$. Note that
\begin{align*}
(\Delta_1 - \sigma \Delta_0)^{-1} \Delta_0 & = (I - (\sigma - \sigma_0)\Delta)^{-1}\Delta = p(\Delta)
\end{align*}
for some polynomial $p$ of degree less than or equal to $2n^2 - n$ because of the Cayley-Hamilton theorem \cite[Theorem 2.4.3.2, Corollary 2.4.3.4]{HJ94_Matrix_Analysis}, and as a result $z_{k+1} = (\Delta_1 - \sigma \Delta_0)^{-1} \Delta_0 = p(\Delta) z_k \in \mathcal{Z}$.

The second approach is inspired by \cite[Theorem 4.7]{CJMU22_linearization} and relies on the property that the invariant vector space $\mathcal{Z}$ is spanned by a subset of the eigenvectors of the pencil $(\Delta_1, \Delta_0)$. More precisely, as shown in  \cref{thm:4_basis}, it is stated that if condition (2) holds, then $z_{k} \in \mathcal{Z}$ can be decomposed as
\begin{align}
\notag
z_{k} = \sum_{i=1}^{n^2 + \ell} \alpha_i \hat{z}_i,
\end{align}
where $\hat{z}_i$ is an eigenvector of \eqref{eq:3_GEP_system} that lies in $\mathcal{Z}$ and $\ell = \frac{1}{2}n(n-1)$. Now $z_{k+1}$ equals,
\begin{align}
\notag
z_{k+1} & = (\Delta_1 - \sigma \Delta_0)^{-1}\Delta_0 z_{k} = \sum_{i=1}^{n^2 + \ell} \frac{\alpha_i}{\lambda_i - \sigma} \hat{z}_i,
\end{align}
and thus $z_{k+1} \in \mathcal{Z}$.
\end{proof}

The consequence of this last statement is that when we choose the initial vector $z_0$ in this invariant subspace $\mathcal{Z}$, then the next iteration will lie in this space as well. As a result, the $\ell$ spurious solutions corresponding to the problem $y^\herm R^\herm M(\lambda, \mu) = 0$ will be filtered. This leads to \cref{alg:Filtering_Arnoldi}:
\begin{algorithm}
\caption{Arnoldi method}
\label{alg:Filtering_Arnoldi}
\begin{algorithmic}[1]
\REQUIRE $z_0 \in \mathcal{Z}$, $\sigma$, $k_{max}$
\STATE{$z_0 \leftarrow z_0 / \|z_0\|_2 $}
\STATE{$Z_1 = \begin{bmatrix}
z_0
\end{bmatrix}$}
\STATE{$H_0 = \begin{bmatrix}
\,
\end{bmatrix}$}
\FOR{$k = 1, 2, \dots, k_{max}$}
\STATE{$\hat{z}_k = (\Delta_1 - \sigma \Delta_0)^{-1}\Delta_0 z_{k-1}$ \hspace{1cm} (Solve using Sylvester equation)}
\STATE{$h_k = Z_{k-1}^\herm \hat{z}_k$}
\STATE{$\hat{z}_k \leftarrow \hat{z}_k - Z_{k-1}h_k$}
\STATE{$\beta_k = \|\hat{z}_k\|_2$}
\STATE{$\hat{z}_k \leftarrow \hat{z}_k/\beta_k$}
\STATE{$\begin{bmatrix}
W_k \\ \hat{V}_k
\end{bmatrix} = \vect^{-1} (\hat{z}_k)$; \quad $V_k = \frac{1}{2}(\hat{V}_k + \hat{V}_k^\tran) \in \mathbb{C}^{n \times n}$; \quad $z_k = \vect \left( \begin{bmatrix}
W_k \\ V_k
\end{bmatrix} \right)$}
\STATE{$Z_{k+1} = \begin{bmatrix}
Z_{k-1} & z_k
\end{bmatrix}$}
\STATE{$H_k = \begin{bmatrix}
H_{k-1} & h_k \\ 0 & \beta_k
\end{bmatrix}$}
\ENDFOR
\RETURN $H_{k_{max}}$, $Z_{k_{max}}$
\end{algorithmic}
\end{algorithm}
Line 5 solves a large linear system of dimension $2n^2 - n$ which is the most expensive step. Because of the structure in the matrices $\Delta_i$, this cost can be significantly reduced by solving the following corresponding generalized Sylvester equation instead:
\begin{align}
\hat{C}X_{k}(A - \sigma B)^\tran - (\hat{A} - \sigma \hat{B}) X_{k} C^\tran = E 
\notag
\end{align}
with $\vect(E) = \Delta_0 z_{k-1}$ and where the solution is $z_{k} = \vect(X_{k})$. If this linear matrix equation is solved using a Bartels-Stewart like method \cite{GLAM92_general_sylvester_method, S16_Review_linear_matrix_equations}, for instance, then the time complexity can be reduced from $\mathcal{O}(n^6)$ to $\mathcal{O}(n^3)$. Lines 6 to 9 in Algorithm \ref{alg:Filtering_Arnoldi} correspond to the orthogonalization step which can be implemented using modified Gram-Schmidt with reorthogonalization to minimize rounding errors. In line $10$, the new vector $z_{k}$ is projected onto $\mathcal{Z}$ to correct for potential deviations caused by floating-point precision errors.

\subsubsection{Two-sided projection}

Instead of projecting the problem to a Hessenberg matrix \eqref{eq:4_Hessenberg_proj}, we also consider the alternative projection
\begin{align}
\notag
\underbrace{Z_k^\herm \Delta_1 Z_k}_{H_1} y = \lambda \underbrace{Z_k^\herm \Delta_0 Z_k}_{H_0} y,
\end{align}
after building the same Krylov space using \cref{alg:Filtering_Arnoldi}, which leads to \cref{alg:Two-sided_Arnoldi}. This new method does not differ much except now the Ritz values are obtained by solving the pencil $H_1 - \lambda H_0$. Both $H_1$ and $H_0$ are Hermitian like the original pencil, so the projection preserves this structure. Moreover, for real eigenvalues $\lambda$ of $\Delta_1 - \lambda \Delta_0$ we expect faster convergence as the column space of $Z_k$ lies closer to both its right and left eigenvector after each iteration. 

Also, an improved filtering of spurious eigenvalues is observed in the numerical experiments, in the sense that twice as many false eigenvalues will be ignored in comparison with the previous filtering. We can explain this using the set of vectors
\begin{align}
\mathcal{W} = \left \{ z \in \mathbb{C}^{2n^2 - n} \mid z = \vect \left( \begin{bmatrix}
W \\ 0
\end{bmatrix} \right) , \, RW = (RW)^\tran \right \} \subset \mathcal{Z}
\label{eq:4_singular_vector_space}
\end{align}
and the following theorem.

\begin{theorem}
\label{prop:4_singular_proj}
If $z \in \mathcal{W}$ \eqref{eq:4_singular_vector_space} and $\col(Z) \subseteq \mathcal{Z}$ \eqref{eq:3_vector_space}, then $Z^\herm (\Delta_1 - \lambda \Delta_0)z = 0$ for any $\lambda \in \mathbb{C}$. Consequence: If such a vector $z \in \col(Z)$, then $(H_1, H_0)$ is a singular pencil where $H_j = Z^\herm \Delta_j Z$ for $j = 1, 2$.
\end{theorem}

\begin{proof}
Define $\Gamma := A - \lambda B$, then
\begin{align}
\notag
\Delta(\lambda) := \Delta_1 - \lambda \Delta_0 = \Gamma \otimes \begin{bmatrix}
0 & R^\herm C \\ CR & Q
\end{bmatrix} - C \otimes \begin{bmatrix}
0 & R^\herm \Gamma \\ \Gamma R & P
\end{bmatrix}.
\end{align}
The matrix vector product $\Delta(\lambda) z$ with $z \in \mathcal{W}$ can be expressed using the $\vect$ operator:
\begin{align}
\notag
\Delta(\lambda) z = \vect \left( \begin{bmatrix}
0 \\ \hat{\Gamma}
\end{bmatrix} \right),
\end{align}
where
\begin{align}
\notag
\hat{\Gamma} = CRW\Gamma^\tran - \Gamma RWC^\tran.
\end{align}
Since $RW$ is symmetric, $\hat{\Gamma}$ is skew symmetric, i.e., $\hat{\Gamma}^\tran = -\hat{\Gamma}$. Let $z_i$ be the $i$-th column of $Z$, then $z_i \in \mathcal{Z}$ because $\col(Z) \subseteq \mathcal{Z}$. Therefore there exists a matrix $W_i \in \mathbb{C}^{(n-1) \times n}$ and a symmetric matrix $V_i = V_i ^\tran \in \mathbb{C}^{n \times n}$ such that
\begin{align}
\notag
z_i = \vect \left( \begin{bmatrix}
W_i \\ V_i
\end{bmatrix} \right).
\end{align}
The inner product of two vectorized matrices can be expressed using the trace operator
\begin{align}
\notag
\alpha_i := z_i^\herm \Delta(\lambda) z = \trace(V_i^\herm \hat{\Gamma}).
\end{align}
Note that since $V_i$ is symmetric, and $\hat{\Gamma}$ is skew symmetric, we have
\begin{align}
\notag
\alpha_i & = -\trace(V_i^\herm \hat{\Gamma}^\tran) ,
& 
(\hat{\Gamma} = -\hat{\Gamma}^\tran)
\\ \notag
& = -\trace(\hat{\Gamma} V_i^\herm) ,
& 
(V_i = V_i^\tran)
\\ \notag
& = -\trace(V_i^\herm \hat{\Gamma}) = -\alpha_i ,
& (\trace(XY) = \trace(YX))
\end{align}
and thus $\alpha_i = 0$. This is true for any column of $Z$ and therefore
\begin{align}
\notag
Z^\herm \Delta (\lambda) z = Z^\herm (\Delta_1 - \lambda \Delta_0) z = 0
\end{align}
for any $\lambda \in \mathbb{C}$.
\end{proof}

A consequence of this theorem is that the eigenvalues of the pencil $(\Delta_1, \Delta_0)$ that have a corresponding eigenvector $z$ in the set $\mathcal{W} \subset \mathcal{Z}$ are no eigenvalues of the projected pencil $(H_1, H_0)$, because in that case $(H_1 - \sigma H_0)y = 0$ for any scalar $\sigma$, with $z = Zy$. Coming back to the characterization of spurious eigenvalues in \Cref{sec:Linearization_Spurious}, note that the solutions $(\lambda, \mu, x)$ of the rMEP $M(\lambda, \mu)Rx$ correspond to a solution $(\lambda, z)$ of the pencil $(\Delta_1, \Delta_0)$ with $z \in \mathcal{W}$, so these $\ell = \frac{1}{2}n(n-1)$ spurious solutions will be filtered by \cref{alg:Two-sided_Arnoldi}. By symmetry, a similar reasoning reveals that the $\ell$ solutions $(\lambda, \mu, x)$ of the pencil $x^\herm R^\herm M(\lambda, \mu)$ are filtered as well and thus we expect $2\ell$ spurious solutions to be ignored.  

The downside of \cref{alg:Two-sided_Arnoldi} is the extra computation required for the projection and the fact that $(H_1, H_0)$ gets closer to a singular pencil after every additional iteration. The issue is that as $Z_k$ grows in size, its column space gets closer to vectors in the set $\mathcal{W}$ that make the pencil singular. In fact, from a certain iteration on, the column space of $Z_k$ in \cref{alg:Two-sided_Arnoldi} will contain a vector $z \in \mathcal{W}$ such that the projected pencil is singular.

\begin{proposition}
Let $Z_k \in \mathbb{C}^{(2n^2 - n) \times k}$ be generated by $k$ iterations of \cref{alg:Filtering_Arnoldi}, $H_i = Z_k^\herm \Delta_j Z_k$ for $j = 0, 1$ and $\col(Z_k) \subseteq \mathcal{Z}$. If $k > n^2$, then the pencil $H_1 - \lambda H_0$ is singular.
\end{proposition}

\begin{proof}
Let $z_i$ be the $i$-th column of $Z_k$, then
\begin{align}
\notag
z_i = \vect \left( \begin{bmatrix}
W_i \\ V_i
\end{bmatrix} \right),
\end{align}
with $V_i = V_i^\tran$ since $\col(Z_k) \subseteq \mathcal{Z}$. From \cref{prop:4_singular_proj}, the pencil is singular if there exists a vector $y = \begin{bmatrix}
y_1 & \dots & y_k
\end{bmatrix}^\tran \in \mathbb{C}^k$ such that
\begin{align}
\notag
Z_ky = \vect \left( \begin{bmatrix}
W \\ 0
\end{bmatrix} \right),
\end{align}
where $RW$ is a symmetric matrix. This is possible if
\begin{align}
\notag
\begin{cases}
\sum_{i=1}^k \Lambda_i y_i = 0, \\
\sum_{i=1}^k V_i y_i = 0,
\end{cases}
\end{align}
with $\Lambda_i := RW_i - W_i^\tran R^\tran$ a skew symmetric matrix. The first equation accounts to $\frac{1}{2}n(n-1)$ equations due to the skew symmetry of $\Lambda_i$, while the second one has $\frac{1}{2}n(n+1)$ equations because of the symmetric $V_i$. Hence, there are a total of $n^2$ linear equations in $k$ unknowns which is guaranteed to have a nonzero solution if $k > n^2$. As a result, the pencil $(H_1, H_0)$ has to be singular for these values of $k$ because of \cref{prop:4_singular_proj}.
\end{proof}

Of course, $k$ does not need to be greater than $n^2$ for the pencil $(H_1, H_0)$ to become singular, as the system we have to solve in the proof of the last corollary could be rank deficient, for example if $Z \subseteq \mathcal{Z}$ contains a vector $z = \vect\left( \begin{bmatrix}
W^\tran & 0
\end{bmatrix}^\tran \right)$ with $RW$ symmetric. It is unlikely for this to happen, but if our Krylov space $Z$ gets closer to one of these vectors $z$, the projected eigenvalue problem $(H_1, H_0)$ becomes close to singular rather quickly. We will illustrate this aspect in \Cref{sec:Experiments}.

\begin{algorithm}
\caption{Two-sided projection method}
\label{alg:Two-sided_Arnoldi}
\begin{algorithmic}[1]
\REQUIRE $z_0 \in \mathcal{Z}$, $\sigma$, $k_{max}$
\STATE{Obtain $Z_{k_{max}}$ from \cref{alg:Filtering_Arnoldi}}
\STATE{Compute $H_j = Z_{k_{max}}^\herm \Delta_j Z_{k_{max}}$ for $j = 0, 1$}
\STATE{Solve for eigenpairs $(\theta_i, x_i)$: $H_1x_i = \theta_i H_0 x_i$}
\RETURN $\theta_i$
\end{algorithmic}
\end{algorithm}

\subsection{Singular pencil}
\label{sec:Numerical_methods_singular}

The methods developed in \Cref{sec:Numerical_methods_nonsingular} assume that $(\Delta_1, \Delta_0)$ is a regular pencil which holds under the conditions of \cref{thm:singular_Delta_0}. In case $\rank(C) < n - 1$ these conditions are not satisfied and, as discussed in the last paragraph of \Cref{sec:Linearization_Compact}, it can be shown that the whole pencil $(\Delta_1, \Delta_0)$ becomes singular. For such an eigenvalue problem, an eigenvalue $\lambda$ is defined as the value for which the rank of $\Delta_1 - \lambda \Delta_0$ is smaller than the normal rank
\begin{align}
\nrank(\Delta_1 - \nu \Delta_0) = \max_{\nu \in \mathbb{C}} \rank( \Delta_1 - \nu \Delta_0).
\notag
\end{align}
Fortunately, the eigenvalues of this GEP corresponding to solutions of the NEPv \eqref{eq:NEPv_QF} hold to this definition under the assumptions of \cref{lemma:4_nrank} and \cref{thm:4_singular_eigenvalue}.

\begin{lemma}
Let $\rank(C) < n-1$ and assume that the eigenvectors of the pencil $(C, B)$ corresponding to nonzero eigenvalues do not lie in the column space of $R$. Let $z$ be a nonzero vector such that $\Delta_0 z = 0$, then $(\Delta_1 - \sigma \Delta_0)z = 0$ for any $\sigma \in \mathbb{C}$.
\label{lemma:4_nrank}
\end{lemma}

\begin{proof}
Let $z = \vect(Z)$ be a nonzero vector such that $\Delta_0 z = 0$ or equivalently
\begin{align}
\hat{C}ZB^\tran - \hat{B}Z C^\tran = 0,
\label{eq:4_lemma_sylvester}
\end{align}
and denote with $C^\tran X = B^\tran X \Theta$ an eigenvalue decomposition of $(C^\tran, B^\tran)$ with $X \in \mathbb{C}^{n \times n}$ and $\Theta = \diag(\theta_1, \dots, \theta_n) \in \mathbb{R}^{n \times n}$, then
\begin{align*}
\hat{C}\hat{Z} - \hat{B}\hat{Z} \Theta = 0,
\end{align*}
where $\hat{Z} = Z B^\tran X$. The nonzero columns $\hat{z}_i$ of $\hat{Z}$ are eigenvectors of the pencil $(\hat{C}, \hat{B})$ corresponding to real-valued eigenvalues $\theta_i$, meaning there must exist nonzero $x_i \in \mathbb{C}^{n-1}$ such that $CRx_i = \theta_i BRx_i$ because of \cref{lemma:3_ChBh_pencil}. From the assumption that no eigenvectors corresponding to nonzero eigenvalues of $(C, B)$ lie in the column space of $R$, we must deduce that $\theta_i = 0$ or $\hat{z}_i = 0$. As a result, $\hat{C}\hat{Z} = 0$ and since $B$ is nonsingular we have $\hat{C}Z = 0$. This implies the second term in \eqref{eq:4_lemma_sylvester} must be zero as well, that is, $\hat{B}ZC^\tran = 0$. Suppose now that $ZC^\tran \ne 0$, then $(\hat{C} - \theta \hat{B})ZC^\tran = 0$ for any value $\theta \in \mathbb{C}$, which is not possible according to \cref{lemma:3_ChBh_pencil} as $(C, B)$ has a finite number of eigenvalues and $R$ is of full rank. Therefore 
\begin{align*}
\hat{C}Z = 0 \quad \text{and} \quad ZC^\tran = 0.   
\end{align*}
It is now easy to see that $\Delta_1 z$ must be zero as well and thus $(\Delta_1 - \sigma \Delta_0)z = 0$ for any shift $\sigma$.
\end{proof}

A direct consequence of this Lemma is that the rank of $(\Delta_1 - \sigma \Delta_0)$ is smaller than or equal to the rank of $\rank(\Delta_0)$, which implies that $\nrank(\Delta_1 - \nu \Delta_0) = \rank(\Delta_0)$. Moreover, both matrices have a common null space which is an important fact that will make the Sylvester equations in \cref{alg:Filtering_Arnoldi} solvable.

Also note that if the nonzero eigenvalues of the pencil $(C, B)$ are simple, the assumption of \cref{lemma:4_nrank} holds for almost all choices of $R$.

\begin{theorem}
Suppose $\rank(C) < n-1$ and that the assumption of \cref{lemma:4_nrank} holds. If $(\lambda, v)$ solves the NEPv \eqref{eq:NEPv_QF} with $v\in \mathbb{C}^n$ such that $Cv$ is nonzero, then $\rank(\Delta_1 - \lambda \Delta_0) < \nrank(\Delta_1 - \sigma \Delta_0)$.
\label{thm:4_singular_eigenvalue}
\end{theorem}

\begin{proof}
Let $(\lambda, v)$ be a solution of the NEPv \eqref{eq:NEPv_QF} then from \cref{prop:2_NEPv_to_NEPv_mu} and \cref{prop:NEPv2MEP} there exists a scalar $\mu \in \mathbb{C}$ and a nonzero vector $w \in \mathbb{C}^{2n-1}\backslash\{0\}$ such that $(\lambda, \mu, v, w)$ solves the MEP \eqref{eq:3_MEP_R}. Note that even though $\Delta_0$ is singular, the eigenpair $(\lambda, z)$ where $z = v \otimes w$ still solves $\Delta_1z = \lambda \Delta_0z$:
\begin{align*}
(\Delta_1 - \lambda \Delta_0) z & = (A - \lambda B) v \otimes \hat{C} w - Cv \otimes (\hat{A} - \lambda \hat{B})w, \\
& = (\mu C)v \otimes \hat{C} w - Cv \otimes (\mu \hat{C})w = 0.
\end{align*}
Suppose this eigenvector also lies in the right null space of $\Delta_0$, then
\begin{align*}
\Delta_0 z = 0,
\end{align*}
or written in terms of Kronecker products
\begin{align*}
B v \otimes \hat{C}w - C v \otimes \hat{B}w = 0.
\end{align*}
Therefore, there must exist a value $\theta \in \mathbb{C}$ such that
\begin{align*}
\begin{cases}
Cv = \theta Bv, \\
\hat{C}w = \theta \hat{B} w.
\end{cases}
\end{align*}
The first equation implies that $\theta$ is real-valued and thus \cref{lemma:3_ChBh_pencil} can be used once more with the second equation, i.e., there exists some $x \in \mathbb{C}^{n-1}$ such that $CRx = \theta BRx$. From the assumption of \cref{lemma:4_nrank}, $\theta$ can only be equal to zero which implies $Cv = 0$. This contradicts with the assumption that $Cv$ is nonzero, hence
\begin{align*}
\Delta_0 z \ne 0.
\end{align*}
Indeed, we have that $z$ lies in the null space of $\Delta_1 - \lambda \Delta_0$ but not in the null space of $\Delta_0$ whereas every vector $\hat{z}$ in the null space of $\Delta_0$ also satisfies $(\Delta_1 - \lambda \Delta_0 ) \hat{z} = 0$ according to \cref{lemma:4_nrank}, and therefore 
\begin{align*}
\rank(\Delta_1 - \lambda \Delta_0) < \rank(\Delta_0) = \nrank (\Delta_1 - \sigma \Delta_0).
\end{align*}
\end{proof}

This theorem indicates that an eigenvalue of the NEPv \eqref{eq:NEPv_QF} is also an eigenvalue of the singular pencil $(\Delta_1, \Delta_0)$ under the two assumptions stated in $\cref{thm:4_singular_eigenvalue}$, allowing us to use dedicated algorithms for solving singular matrix pencils \cite{HMP23_Singular_GEP_part_II, MW24_shift_and_invert_Arnoldi}. An observation is that the standard Arnoldi algorithm (\cref{alg:Filtering_Arnoldi} without the projection onto $\mathcal{Z}$ in line 10) can be used as well: The system 
\begin{align}
(\Delta_1 - \sigma \Delta_0) \hat{z}_k = \Delta_0 z_{k-1}
\label{eq:4_sing_sys}
\end{align}
that needs to be solved in each iteration is singular, but if the shift $\sigma$ is not an eigenvalue then both $\Delta_1 - \sigma \Delta_0$ and $\Delta_0$ have the same column space as a result of \cref{lemma:4_nrank} and the fact that both $\Delta_1$ and $\Delta_0$ are Hermitian\footnote{If $\sigma$ is not an eigenvalue then neither is the complex conjugate of $\sigma$ because both $\Delta_0$ and $\Delta_1$ are Hermitian. The null spaces of $\Delta_1 - \overline{\sigma} \Delta_0$ and $\Delta_0$ must be equal as a consequence of \cref{lemma:4_nrank}, and thus their row spaces must be equal as well. As a result, $\Delta_1 - \sigma \Delta_0$ and $\Delta_0$ share the same column space.}. In other words, the Sylvester equation equivalent to \eqref{eq:4_sing_sys} has an infinite number of solutions of which we can pick one to continue the Arnoldi iterations. Unfortunately, the filtering proposed in \Cref{sec:Numerical_methods_nonsingular} is not possible due to the fact that the singular pencil $(\Delta_1, \Delta_0)$ will not have enough eigenvectors to form a basis for the invariant subspace $\mathcal{Z}$ \eqref{eq:3_vector_space}. Moreover, the upper bound on the number of solutions of the NEPv \eqref{eq:NEPv_QF} drops from $n^2$ to $n(2r+1) - r(r+1)$ with $r = \rank(C)$ (\cref{thm:bkk}), which can be seen as a disadvantage as the difference between the number of NEPv solutions and the dimension of $\Delta_0$ may become larger. A proof for \cref{thm:bkk} is given in \Cref{app:upper_bound}.

\begin{theorem}
Let $r = \rank(C) < n-1$, then the number of isolated solutions $(\lambda, \mu) \in \mathbb{C}\backslash\{0\} \times \mathbb{C}\backslash\{0\}$ of the system of polynomial equations \eqref{eq:NEPv_polyn}
\begin{align}
\begin{cases}
f(\lambda, \mu) = \det(A - \lambda B - \mu C) = 0,
\\
g(\lambda, \mu) = \trace \big( (P - \mu Q) \, \adj(A - \lambda B - \mu C) \big) = 0,
\end{cases}
\notag
\end{align}
is bounded by $n(2r + 1) - r(r + 1)$.
\label{thm:bkk}
\end{theorem}

\section{Numerical examples}
\label{sec:Experiments}

This section is concerned with solving the NEPv \eqref{eq:NEPv_QF} for different matrices using the two algorithms from \Cref{sec:Numerical_methods} to solve the exact linearization in \Cref{sec:Linearization}. Recall that the linearization is, 
\begin{align}
\Delta_1 z = \lambda \Delta_0 z,
\notag
\end{align}
with,
\begin{gather}
\Delta_1 = A \otimes \hat{C} - C \otimes \hat{A}, \qquad \Delta_0 = B \otimes \hat{C} - C \otimes \hat{B},
\notag
\\
\hat{A} = \begin{bmatrix}
0 & R^\herm A \\ AR & P
\end{bmatrix}, \quad \hat{B} = \begin{bmatrix}
0 & R^\herm B \\ BR & 0
\end{bmatrix}, \quad \hat{C} = \begin{bmatrix}
0 & R^\herm C \\ CR & Q
\end{bmatrix}
\notag
\end{gather}
where the matrix $R \in \mathbb{C}^{n \times (n-1)}$ can be chosen at random, but of full rank. The NEPv matrices $A, B, C, P$ and $Q$ are symmetric matrices of size $n$ and additionally $B$ and $Q$ are positive definite. The example code used to generate the figures is available at \url{https://gitlab.kuleuven.be/numa/public/nepv}.


\subsection{Example 1}
\label{sec:ex1}

Consider the NEPv of dimension $n = 100$ where the matrices are dense and randomly generated as 
\begin{lstlisting}[style=Matlab-editor]
rng(0)
A = randn(n) + 1i*randn(n); A = (A + A')/2;
B = randn(n) + 1i*randn(n); B = sqrtm(B*B');
C = randn(n) + 1i*randn(n); C = (C + C')/2;
P = randn(n) + 1i*randn(n); P = (P + P')/2;
Q = randn(n) + 1i*randn(n); Q = sqrtm(Q*Q');
\end{lstlisting}
in \textsc{Matlab}. We choose random complex numbers as the entries for the matrix $R \in \mathbb{C}^{n \times (n-1)}$, and compute the matrices $\hat{A}$, $\hat{B}$ and $\hat{C}$ of dimension $2n-1 = 199$, resulting in operator determinants $\Delta_i$ of size $2n^2 - n = 19900$. These delta matrices do not need to be computed explicitly since the linear system in \cref{alg:Filtering_Arnoldi} is solved as a Sylvester equation using the MEP matrices. For this example we set $\sigma$ equal to zero to search for the eigenvalues with smallest magnitude, and we take a random starting vector $z_0$ in $\mathcal{Z}$.

\begin{figure}
\centering
\includegraphics[width=0.8\linewidth]{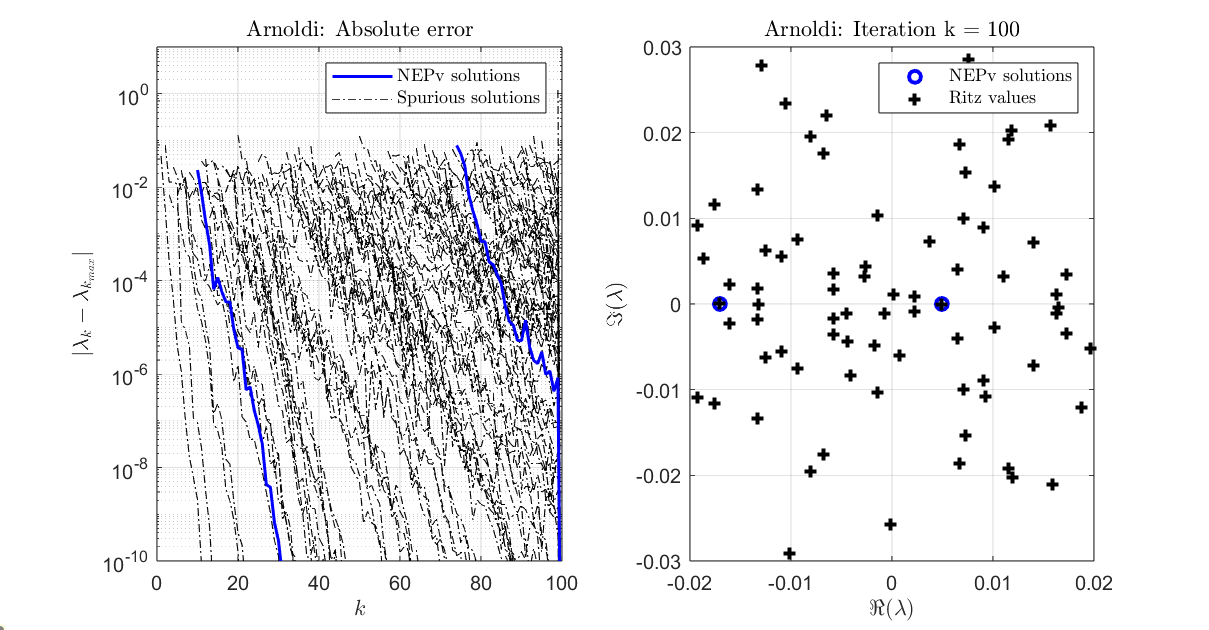}
\caption{Left figure: Convergence plot for the example in \Cref{sec:ex1}. The blue lines represent convergence of Ritz values computed in \cref{alg:Filtering_Arnoldi} towards eigenvalues of the NEPv, while the other curves represent convergence of Ritz values towards spurious solutions. Right figure: Comparison between the Ritz values and the NEPv eigenvalues at iteration $k=100$.}
\label{fig:n_100_Arnoldi}
\end{figure}

When the problem is solved using \cref{alg:Filtering_Arnoldi}, we get convergence to the eigenvalues as shown in the left figure of \cref{fig:n_100_Arnoldi} where the absolute error between the Ritz values of the current iteration and the reference eigenvalues are plotted in function of the iteration $k$. The reference eigenvalues are chosen as the Ritz values at the last iteration $k_{max}$, which explains the sudden vertical jump in the figures at $k = 100$. The blue convergence curves are of interest as they correspond to  Ritz values that approximate an eigenvalue of the NEPv. Convergence to the first eigenvalue happens around iteration $k=30$, while the second blue line has not converged yet in iteration $k=100$. The right plot in \cref{fig:n_100_Arnoldi} compares the computed Ritz values in the last iteration with the eigenvalues of the NEPv.

\begin{figure}
\centering
\includegraphics[width=0.8\linewidth]{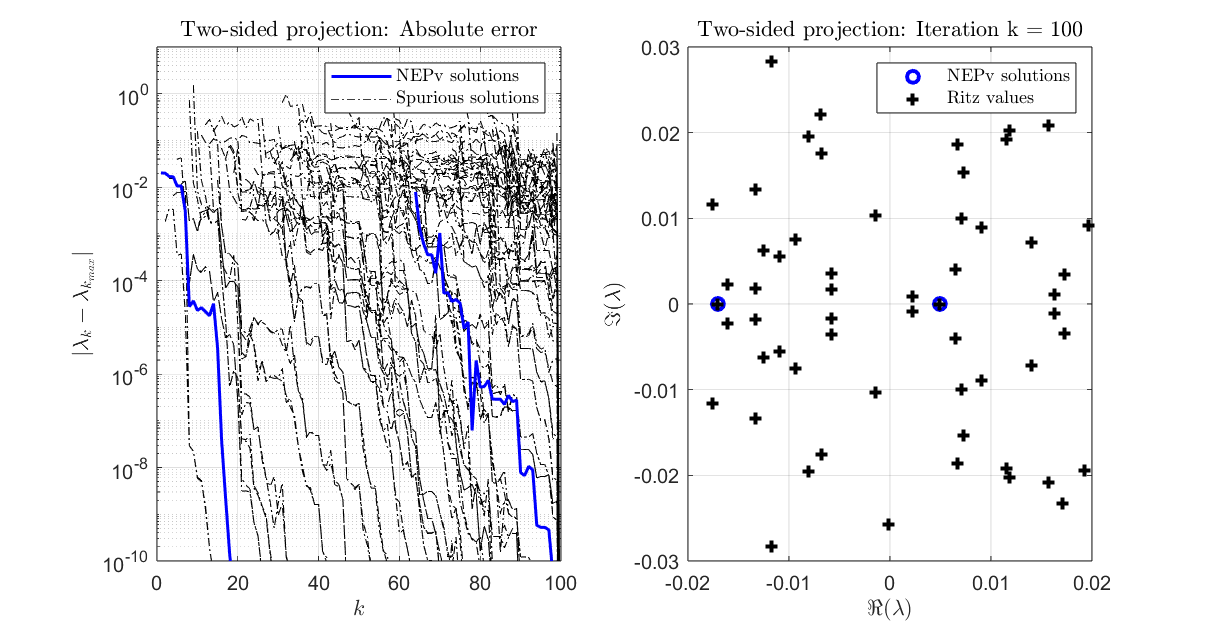}
\caption{Similar to \cref{fig:n_100_Arnoldi} except now \cref{alg:Two-sided_Arnoldi} is used to compute the Ritz values.}
\label{fig:n_100_two_sided}
\end{figure}

\begin{figure}
\centering
\includegraphics[width=0.5\linewidth]{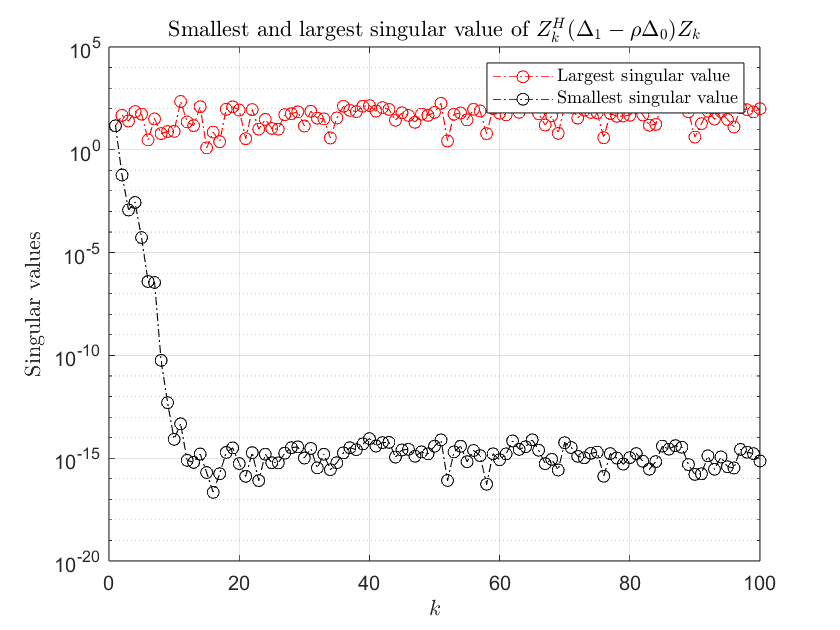}
\caption{Illustrates after which iteration the projected pencil $H_1 - \lambda H_0$ becomes numerically rank deficient for the example in \Cref{sec:ex1}.}
\label{fig:n_100_singularity}
\end{figure}

Alternatively, if \cref{alg:Two-sided_Arnoldi} is used where we also compute the two-sided projection to obtain the pencil $H_1 - \lambda H_0$, then convergence to the real-valued eigenvalues is reached within fewer iterations ($k = 18$ and $k = 97$) as shown in \cref{fig:n_100_two_sided}. Comparing the right plots of \cref{fig:n_100_Arnoldi} and \cref{fig:n_100_two_sided}, we see fewer Ritz values for this second method due to the improved filtering. The drawback of this approach is the extra computation needed to obtain the projected pencil, and the fact that this pencil becomes singular after a certain number of iterations. This is illustrated in \cref{fig:n_100_singularity} by computing the smallest singular values of $H_0 + \rho H_1$ for some randomly chosen value $\rho$. Note that the smallest singular value drops to machine precision around $k = 12$, meaning the projected pencil is numerically rank deficient at this point.

\subsection{Example 2}
\label{sec:ex2}

Consider the scalar function $u(x) \in \mathbb{R}$ in function of $x \in \mathbb{R}$ such that
\begin{align*}
-\nabla^2 u(x) + f(u) c(x) u(x) = \lambda u(x).
\end{align*}
In this example we solve a large sparse system that arises by discretizing this wave equation using finite differences on the domain $[-L/2, L/2]$ with Dirichlet boundary conditions $u(-L/2) = u(L/2) = 0$. The goal is to find the smallest eigenvalues $\lambda$ and their corresponding wave functions $u(x)$ such that $\int_{-L/2}^{L/2} |u(x)|^2 dx = 1$. This equation is similar to the problem in \cite{CJMU22_linearization}, except now the nonlinearity is defined as
\begin{align*}
f(u) = \int_{-L/2}^{L/2} p(x) |\nabla u(x)|^2 dx, \quad \text{with} \quad p(x) = 5 \cos(\frac{\pi}{L}x),
\end{align*}
and we take
\begin{align*}
c(x) = 1 - \exp \left( -\frac{(10x - 1)^2}{10}\right).
\end{align*}
The equation is discretized using central differences and an equidistant mesh of $n+2$ nodes $\{x_k\}_{k=0}^{n+1}$ with $x_0 = -L/2$ and $x_{n+1} = L/2$, where the distance between two nodes is $h = L/(n+1)$. The discrete wavefunction is stored in the vector $v \in \mathbb{C}^n$ such that $v_i$ corresponds to $u(x_i)$ for $i = 1, \dots, n$. The resulting discretized problem is
\begin{align*}
Av = \lambda v + \frac{v^\herm P v}{v^\herm v} Cv,
\end{align*}
with banded matrices computed as
\begin{gather*}
A_{i,i} = 2/h^2, \quad A_{i, i+1} = A_{i+1, i} = -1/h^2, \\
C = - \diag (\begin{bmatrix}
c(x_1) & c(x_2) & \dots & c(x_n)
\end{bmatrix}), \\
P_{i, i} = \frac{p(x_{i-1}) + p(x_{i+1})}{4h^2}, \quad P_{i, i+2} = P_{i+2, i} = \frac{-p(x_{i+1})}{4h^2}.
\end{gather*}

\begin{figure}
\centering
\includegraphics[width=0.8\linewidth]{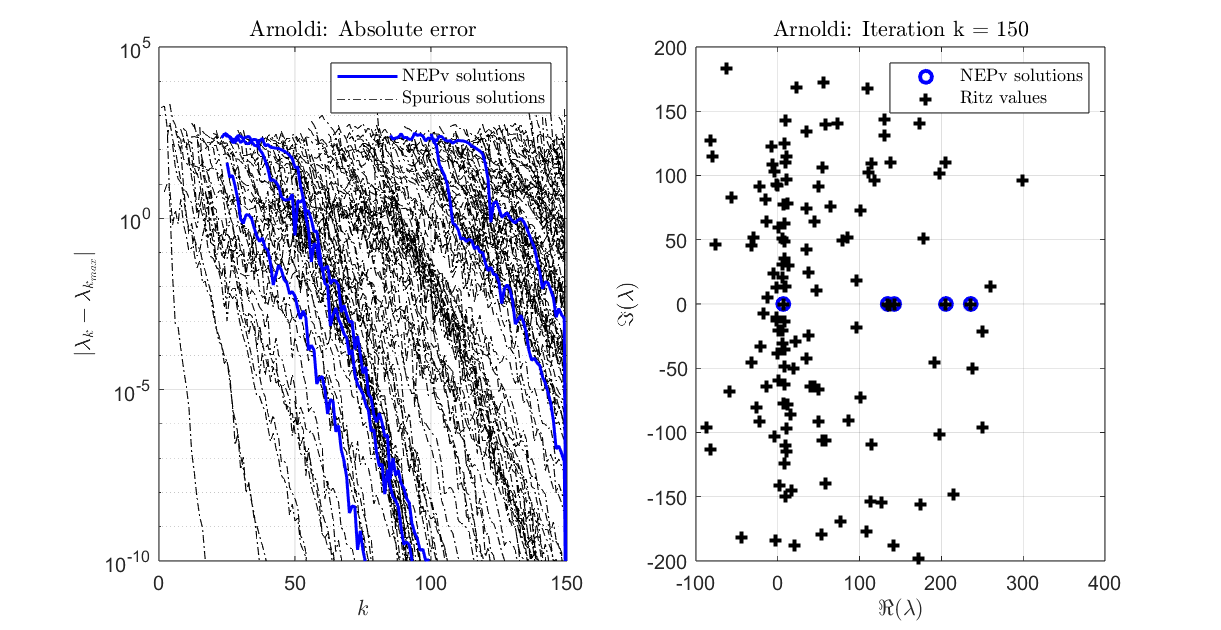}
\caption{Left figure: Convergence plot for the example in \Cref{sec:ex2}. The blue lines represent convergence of Ritz values computed in \cref{alg:Filtering_Arnoldi} towards eigenvalues of the NEPv, while the other curves represent convergence of Ritz values towards spurious solutions. Right figure: Comparison between the Ritz values and the NEPv eigenvalues at iteration $k=150$.}
\label{fig:n_256_Arnoldi}
\end{figure}

\begin{figure}
\centering
\includegraphics[width=0.8\linewidth]{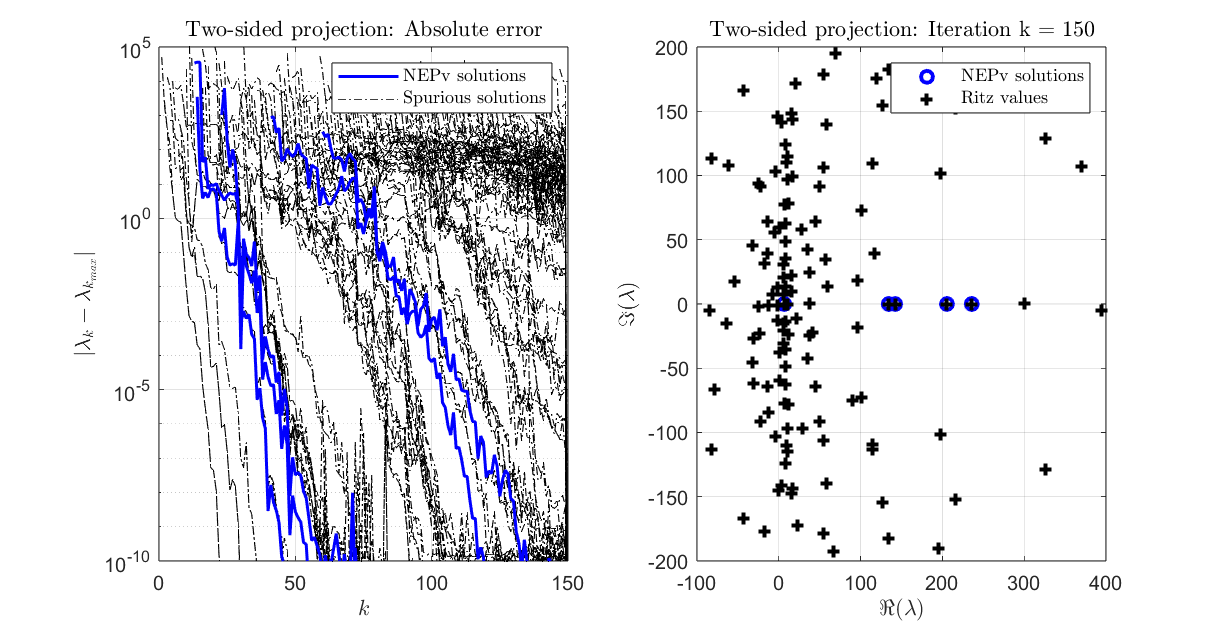}
\caption{Similar to \cref{fig:n_256_Arnoldi} except now \cref{alg:Two-sided_Arnoldi} is used to compute the Ritz values.}
\label{fig:n_256_two-sided}
\end{figure}

\begin{figure}
\centering
\includegraphics[width=0.5\linewidth]{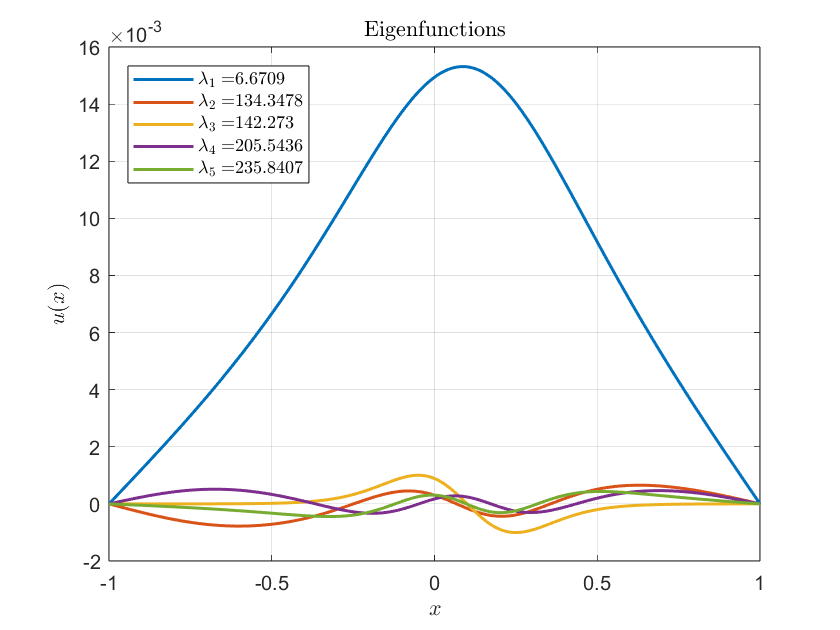}
\caption{Eigenfunctions $u(x)$, scaled with its eigenvalue, corresponding to the 5 computed eigenvalues in function of $x$ for the example in \Cref{sec:ex2}.}
\label{fig:n_256_eigenfunctions}
\end{figure}

For this example, we choose the values $L = 2$ and $n = 256$ and use the shift $\sigma = 50$. The matrix $R$ is chosen as a horizontal stack of the identity matrix and a row vector with random complex elements. After 150 iterations of both \cref{alg:Filtering_Arnoldi} and \cref{alg:Two-sided_Arnoldi}, convergence has been reached to 4 and 5 eigenvalues, respectively. The approximated absolute errors of the computed eigenvalues are plotted in \cref{fig:n_256_Arnoldi} and \cref{fig:n_256_two-sided}.
In this case, \cref{alg:Two-sided_Arnoldi} has converged to the smallest eigenvalue $\lambda = 6.67$ in 45 iterations, showing its faster convergence rate to the NEPv eigenvalues in comparison with \cref{alg:Filtering_Arnoldi} which converges to the same eigenvalue after 75 iterations. The entries of the eigenvectors corresponding to the $5$ computed eigenvalues are the discretization of the eigenfunctions $u(x)$ and are plotted in \cref{fig:n_256_eigenfunctions}.

\subsection{Example 3}
\label{sec:ex3}

\begin{figure}
\centering
\includegraphics[width=0.8\linewidth]{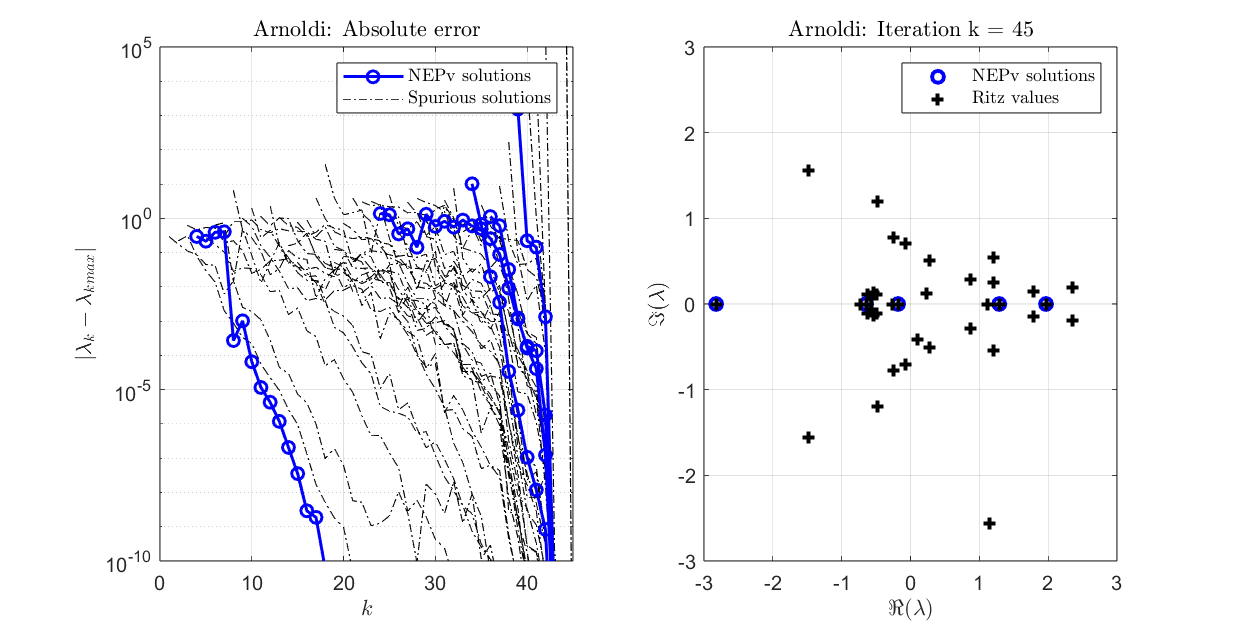}
\caption{Left figure: Convergence plot for the example in \Cref{sec:ex3}. The blue lines represent convergence of Ritz values computed using standard Arnoldi towards eigenvalues of the NEPv, while the other curves represent convergence of Ritz values towards spurious solutions. Right figure: Comparison between the Ritz values and the NEPv eigenvalues at iteration $k=42$.}
\label{fig:n_5_singular_Arnoldi}
\end{figure}

Consider the NEPv \eqref{eq:NEPv_QF} with a rank $r$ matrix $C$ generated in \textsc{Matlab} as
\begin{lstlisting}[style=Matlab-editor]
C1 = orth(randn(n, r) + 1i*randn(n, r));
C = C1*diag(randn(r, 1))*C1';
\end{lstlisting}
and where the other matrices are generated as in \Cref{sec:ex1}. For this example we take $n = 5$, $r = 2$, a shift $\sigma = 0$ and a complex randomly generated matrix $R$. To retrieve the eigenvalues, we solve the pencil $(\Delta_1, \Delta_0)$ using the standard Arnoldi method, i.e., \cref{alg:Filtering_Arnoldi} without the projection step onto the vector space $\mathcal{Z}$ in line 10. The singular system in line 5 is solved using the Bartels-Steward method \cite{GLAM92_general_sylvester_method}, in which the smaller singular systems are solved by using an LU factorization and discarding the zero rows. The convergence behavior for all Ritz values is plotted in \cref{fig:n_5_singular_Arnoldi}. Convergence to the smallest eigenvalue is reached around iteration $k = 18$.

\section{Conclusion and outlook}
\label{sec:Conclusion}

The NEPv considered in this paper contains one scalar nonlinearity in which the eigenvector appears quadratically. We show how the solutions of this problem are related to the real-valued solutions of a system of polynomial equations, and that the upper bound on the number of isolated solutions is $n^2$ where $n$ is the dimension of the NEPv. The first main contribution of this paper is the conversion to a GEP in the sense that the spectrum of the resulting linear problem contains all the eigenvalues of the NEPv. We stated under which conditions the resulting GEP may become singular. As a second contribution, we show how the Arnoldi method can be used to efficiently obtain the eigenvalues of the NEPv by exploiting the structure of the resulting linear problem and by filtering a portion of the spurious solutions. Using the alternative two-sided projection results in an improved filtering of the spurious solutions at the cost of explicitly having to compute the projection. In case the linear pencil is singular because of a low rank matrix $C$, the methods can still be used since system \eqref{eq:4_sing_sys} is solvable whenever the shift is not an eigenvalue, however the filtering approach is not applicable anymore. The theoretical results are supported by numerical examples.

The paper only discussed the NEPv with one scalar nonlinearity, but this can be generalized to multiple terms in future research. A problem that will occur however is the exponential increase in the size of the resulting linear problem as the number of scalar nonlinearities is increased.

\appendix

\section{Technical Lemmas}
\label{sec:app_tech_lemmas}

\begin{lemma}
Assume $\Delta_0$ is nonsingular and that \eqref{eq:3_GEP_system} has $2n^2 - n$ distinct eigentuples $(\lambda, \mu) \in \mathbb{C} \times \mathbb{C}$. Let $(\lambda, \mu)$ and $z \in \mathbb{C}^{2n^2 - n}$ be such an eigentuple with its corresponding eigenvector, then at least one of the following statements must be true:
\begin{itemize}
\item The eigenvector $z$ is an element of $\mathcal{Z}$.
\item There exists a nonzero vector $y \in \mathbb{C}^{n-1}$ such that $y^\herm R^\herm M(\lambda, \mu) = 0$.
\end{itemize}
\label{thm:4_invariant_subspace}
\end{lemma}

\begin{proof}
Because $\Delta_0$ is nonsingular, \cref{thm:MEP2GEP} says there must exist a $v$ and $w$ such that $(\lambda, \mu, v, w)$ solves the MEP \eqref{eq:3_MEP_R} with $z = v \otimes w$, i.e.
\begin{align}
\begin{cases}
Mv = 0, \\
R^\herm M w_2 = 0, \\
MRw_1 + Sw_2 = 0,
\end{cases}
\label{eq:4_MEP_eq}
\end{align}
where $w = \begin{bmatrix}
w_1^\tran & w_2^\tran
\end{bmatrix}^\tran$. Note that if $\rank(M) < n-1$, the eigenpair has a geometric multiplicity greater than one which contradicts with the assumption that all of them are distinct, thus $\rank(M) = n-1$. We consider the following two cases:
\begin{itemize}
\item Case I: $\rank(R^\herm M) < n-1$. The matrix $R^\herm M$ is rank deficient and there must exist a nonzero vector $y$ in its left null space, i.e., $y^\herm R^\herm M = 0$.
\item Case II: $\rank(R^\herm M) = n-1$. In this case, $M$ and $R^\herm M$ share the same one-dimensional (right) null space and from \eqref{eq:4_MEP_eq} both $v$ and $w_2$ lie in this space, therefore there must exist an $\alpha \in \mathbb{C}$ such that $w_2 = \alpha v$. Now we have
\begin{align*}
z & = v \otimes w = v \otimes \begin{bmatrix}
w_1 \\ \alpha v
\end{bmatrix} = \vect \left( \begin{bmatrix}
w_1 v^\tran \\ \alpha v v^\tran
\end{bmatrix} \right) \in \mathcal{Z}.
\end{align*}
\end{itemize}
\end{proof}

\begin{lemma}
Suppose $\Delta_0$ is nonsingular and that \eqref{eq:3_GEP_system} has $2n^2 - n$ distinct eigentuples $(\lambda, \mu) \in \mathbb{C} \times \mathbb{C}$. Denote $\{(\lambda_i, \mu_i, z_i)\}_{i = 1}^N$ as the set of eigenpairs and corresponding eigenvectors of \eqref{eq:3_GEP_system} for which $z_i \in \mathcal{Z}$, then this set is a basis for $\mathcal{Z}$. 
\label{thm:4_basis}
\end{lemma}

\begin{proof}
The system of GEPs \eqref{eq:3_GEP_system} has distinct eigentuples, so the set of $N$ eigenvectors $\{z_i\}_{i=1}^N$ must be linearly independent and $N \leq \dim(\mathcal{Z})$. Denote the set $\{(\hat{\lambda}_i, \hat{\mu}_i, \hat{z}_i)\}_{i = 1}^{2n^2 - n - N}$ as the other solutions of \eqref{eq:3_GEP_system} for which the eigenvector $\hat{z}_i$ does not lie in $\mathcal{Z}$, then according to \cref{thm:4_invariant_subspace} there must exist a set of nonzero vectors $\{y_i\}_{i=1}^{2n^2 - n - N}$ such that $(\hat{\lambda}_i, \hat{\mu}_i, y_i)$ solves the rMEP $y_i^\herm R^\herm M(\hat{\lambda}_i, \hat{\mu}_i) = 0$. From \cref{thm:singular_Delta_0} and the assumption that $\Delta_0$ is nonsingular, we have that $\rank(CR - \sigma BR) = n-1$ for any $\sigma \in \mathbb{C}$ and if this fact is combined with \cite[Lemma 1]{HKP23_rMEP}, we must deduce that the number of solutions of this rMEP precisely equals $\ell = \frac{1}{2}n(n-1)$ counting multiplicities, hence
\begin{align*}
2n^2 - n - N \leq \ell, \\
\iff \quad N \geq n^2 + \ell = \dim(\mathcal{Z}).
\end{align*}
Consequently, $N$ must be equal to $\dim(\mathcal{Z})$ and therefore the set $\{z_i\}_{i=1}^N$ forms a basis for $\mathcal{Z}$.
\end{proof}

\section{Proof of \cref{thm:bkk}} 
\label{app:upper_bound}

The proof consists of three parts: In the first two parts the Newton polytopes of both polynomials are determined. In the third part, the Bernstein-Khovanskii-Kushnirenko (BKK) theorem  \cite[Theorem 5.4, Page 346]{CLO04_Using_algebraic_geometry} is used to obtain an upper bound on the number of isolated solutions.

Define $\Gamma = A - \lambda B$ and let $C = GDG^\herm$ be the reduced eigenvalue decomposition such that $D \in \mathbb{C}^{r \times r}$ contains the $r$ nonzero eigenvalues of $C$, and $G \in \mathbb{C}^{n \times r}$ the corresponding eigenvectors. Assume $\Gamma$ is invertible, then from the matrix determinant lemma (see also \cite[Equation (0.8.5.1)]{HJ94_Matrix_Analysis}) we have,
\begin{align}
f(\lambda, \mu) = \det(D) \det(\Gamma)  \det(D^{-1} - \mu G^\herm \Gamma^{-1} G).
\notag
\end{align}
The variable $\mu$ only appears in the last determinant of this expression and since $G^\herm \Gamma^{-1} G$ is a matrix function of $\lambda$ of size $r$, $f(\lambda, \mu)$ can be seen as a function of degree $r$ in $\mu$ with coefficients in function of $\lambda$. Thus if $\lambda^i \mu^j$ is a monomial of $f$, then $j \leq r$. From \cref{prop:2_number_sols}, $f$ has a maximal degree of $n$ and therefore $i + j \leq n$. As a result, the Newton polytope of $f(\lambda, \mu)$ is $P_f = \Conv((0, 0), (n, 0), (n-r, r), (0, r))$. This remains true in case $\Gamma$ is singular since this occurs for a finite number of values $\lambda$ and $f(\lambda, \mu)$ is continuous.

The second function $g(\lambda, \mu)$ can be seen as a function of degree $n-1$ in $\lambda$ with coefficients in function of $\mu$. Alternatively, using the Woodbury matrix identity and the matrix determinant lemma, the adjugate in the expression of $g$ can be written as
\begin{align}
\adj(\Gamma - \mu C) & = f(\lambda, \mu) \Gamma^{-1} + \mu \det(D) \det(\Gamma) \Gamma^{-1} G \adj (D^{-1} - \mu G^\herm \Gamma^{-1} G) G^\herm \Gamma^{-1},
\notag
\end{align}
assuming both $\Gamma$ and $M = \Gamma - \mu C$ are nonsingular, but because this expression is continuous in $\mu$ it also holds in case $M$ is singular. The expression can be seen as a matrix polynomial in $\mu$ of degree $r$ with coefficients in function of $\lambda$, meaning $g(\lambda, \mu)$ has a degree of $r+1$ in $\mu$. Additionally, from \cref{prop:2_number_sols}, this polynomial has degree $n$, thus if $\lambda^i \mu^j$ is a monomial in $g(\lambda, \mu)$, then $i \leq n-1$, $j \leq r + 1$ and $i+j \leq n$. The corresponding Newton polytope is $P_g = \Conv((0, 0), (n-1, 0), (n-1, 1), (n-r-1, r+1), (0, r+1))$. This also remains true in case $\Gamma$ is singular due to the continuity of $g(\lambda, \mu)$.

The BKK theorem \cite[Theorem 5.4, Page 346]{CLO04_Using_algebraic_geometry} states that the number of isolated solutions of the system of polynomial equations $f(\lambda, \mu) = g(\lambda, \mu) = 0$ is bounded by the mixed volume $MV_2(P_f, P_g)$, where $P_f$ and $P_g$ are the Newton polytopes of $f$ and $g$, respectively. This mixed volume is the coefficient of the monomial $\lambda_f \lambda_g$ in
\begin{align}
\Vol_2(\lambda_f P_f + \lambda_g P_g) & = \left(rn - \frac{r^2}{2} \right) \lambda_f^2 + \left( (r + 1)(n - 1) - \frac{r^2}{2} \right) \lambda_g^2
\notag
\\
& + (n + r(2n - r - 1)) \lambda_f \lambda_g,
\notag
\end{align}
and therefore the number of isolated solutions is bounded by $n + r(2n - r - 1)$.

\section*{Acknowledgments}
The authors would like to thank Tom Kaiser for his expertise and insightful discussions, which contributed to the proof of \cref{thm:bkk}.

\bibliographystyle{siamplain}
\bibliography{references}
\end{document}